\documentclass{article}
\usepackage{amsmath, amsthm, amssymb}
\usepackage[pdftex]{graphicx}
\usepackage{multirow}

\title{Spectrum and Heat Kernel Asymptotics on General Laakso Spaces}
\author{Matthew Begu\'{e}, Levi DeValve, David Miller\\ and Benjamin Steinhurst\footnote{Research supported by NSF grant DMS-0505622}}

\numberwithin{equation}{section}

\newtheorem{theorem}{Theorem}[section]
\newtheorem{cor}{Corollary}[section]
\newtheorem{lemma}{Lemma}[section]
\newtheorem{prop}{Proposition}[section]

\newtheorem{defn}{Definition}[section]

\begin{document}
\maketitle

\begin{center}
\small
Contacts:\\
matthew.begue@uconn.edu\\
Matthew Begu\'{e}\\
Department of Mathematics\\
University of Connecticut\\
Storrs, CT 06269 USA
\vskip 3mm
levi.devalve@uconn.edu\\
Levi DeValve\\
Department of Mathematics\\
University of Connecticut\\
Storrs, CT 06269 USA
\vskip 3mm
david.miller@salve.edu\\
David Miller\\
Department of Mathematics\\
Salve Regina University\\
Newport, RI 02840 USA
\vskip 3mm
steinhurst@math.uconn.edu\\
Benjamin Steinhurst\footnote{Corresponding author}\\
Department of Mathematics\\
University of Connecticut\\
Storrs, CT 06269 USA\\
t. +1 (860) 486-3923\\
f. +1 (860) 486-4238
\end{center}

\begin{abstract}
We introduce a method of constructing a general Laakso space while calculating the spectrum and multiplicities of the Laplacian operator on it.  Using this information, we find the leading term of the trace of the heat kernel and the spectral dimension on an arbitrary Laakso space.  
\end{abstract}

\section{Introduction}

Much work has been done on the analysis of fractals, specifically concentrating on the spectrum of the Laplacian operator on irregular domains. One such topic are drums with Koch snowflake boundary, see for example \cite{Lapidus1995}. This paper will be concerned instead with the irregular domain being a fractal itself.   Some notable works with this type of domain include \cite{spectralSG,mtree,kajino,eigen,ben,tutorial} among others.  Laakso's spaces were introduced in \cite{Laakso}.  They are a family of fractals with an arbitrary Hausdorff dimension greater than one and were considered originally for their nice analytic properties.  Constructions of the Laakso spaces are given in \cite{Laakso,eigen,ben} as well as in Section \ref{sec:Laakso} of this paper.  Theorem 6.1 in \cite{eigen} gives the spectrum of the Laplacian operator on any given Laakso space, in Theorem \ref{thm:spect} we give the multiplicities.

%Rewrite this paragraph
An important part of the analysis of Laplacians is the heat equation and associated heat kernel, which can reveal significant information about the operator and underlying space.  The information gained from studying heat kernels can be applied in other areas of analysis as well as other fields such as physics.  The papers \cite{physics,tutorial} are devoted to finding and analyzing the heat kernel and the trace of the heat kernel.  The notion of complex valued fractal dimensions  and the accompanying oscillating behavior of the heat kernel were studied in \cite{physics, spectralSG}.

We begin by reviewing the construction of the Laakso spaces as presented in \cite{Laakso,eigen, ben} in Section \ref{sec:Laakso}.  This section also contains background information on the Hausdorff dimension, its calculation for Laakso spaces, and some specific values for certain Laakso spaces.  In subsection \ref{sec:Laplacian} we define  the Laplacian operator that will be used throughout the rest of this paper.

In Section \ref{sec:delta} we begin by stating the spectrum of the Laplacian and the multiplicities of each eigenvalue (Theorem \ref{thm:spect}), while the rest of the section is devoted to the proof of this result.   In Section \ref{counts} we provide an analytical proof by examining, as in \cite{eigen}, the different ``shapes" that make up the space.  Since each shape has a unique contribution to the spectrum  counting the number of shapes allows us to calculate the spectrum with multiplicities.  In Section \ref{ssec:matlab} we verify the results computationally using MATLAB.  Finally in Sections \ref{sec:heatmethods} and \ref{sec:heat} we use the spectrum and multiplicities obtained in  Sections \ref{hdim} and \ref{counts} to calculate the trace of the heat kernel using the same method outlined for diamond fractals in \cite{physics}.

\section{Laakso Spaces}\label{sec:Laakso}
The spaces that will be analyzed were first defined by Laakso in \cite{Laakso}.  Laakso's spaces form an uncountable family of metric-measure spaces indexed by sequences $\{j_n\}_{n=1}^\infty$.  An equivalent construction using projective limits was hinted at in \cite{BarlowEvans2004} and fully developed in \cite{ben}.  Then in \cite{eigen}, a more in depth description of the projective limit construction was used to calculate the spectrum of the Laplacian constructed in \cite{ben}.  We will be using the construction in \cite{eigen} as it is also well-suited for our calculations.

The Laakso space can be visualized with a sequence of quantum graphs, denoted $F_n,n\geq0$, each an increasingly better approximation of the Laakso space.  The first of these graphs is the unit interval, denoted $F_0$.  Laakso spaces are defined by a sequence $\{j_n\}_{n=1}^\infty$ of integers $j_n\geq2$, where each $j_n$ described the number of identifications at step $n$ of the construction.  To construct the graph of $F_{n+1}$, first every cell, or interval between two nodes, of the $F_n$ graph is split evenly into $j_n$ segments by adding nodes.  This graph is then duplicated and connected at the newly-added notes.  In this visualization, all nodes are arranged in columns.

\begin{figure}[htbp]
\begin{center}
\includegraphics[scale=.25]{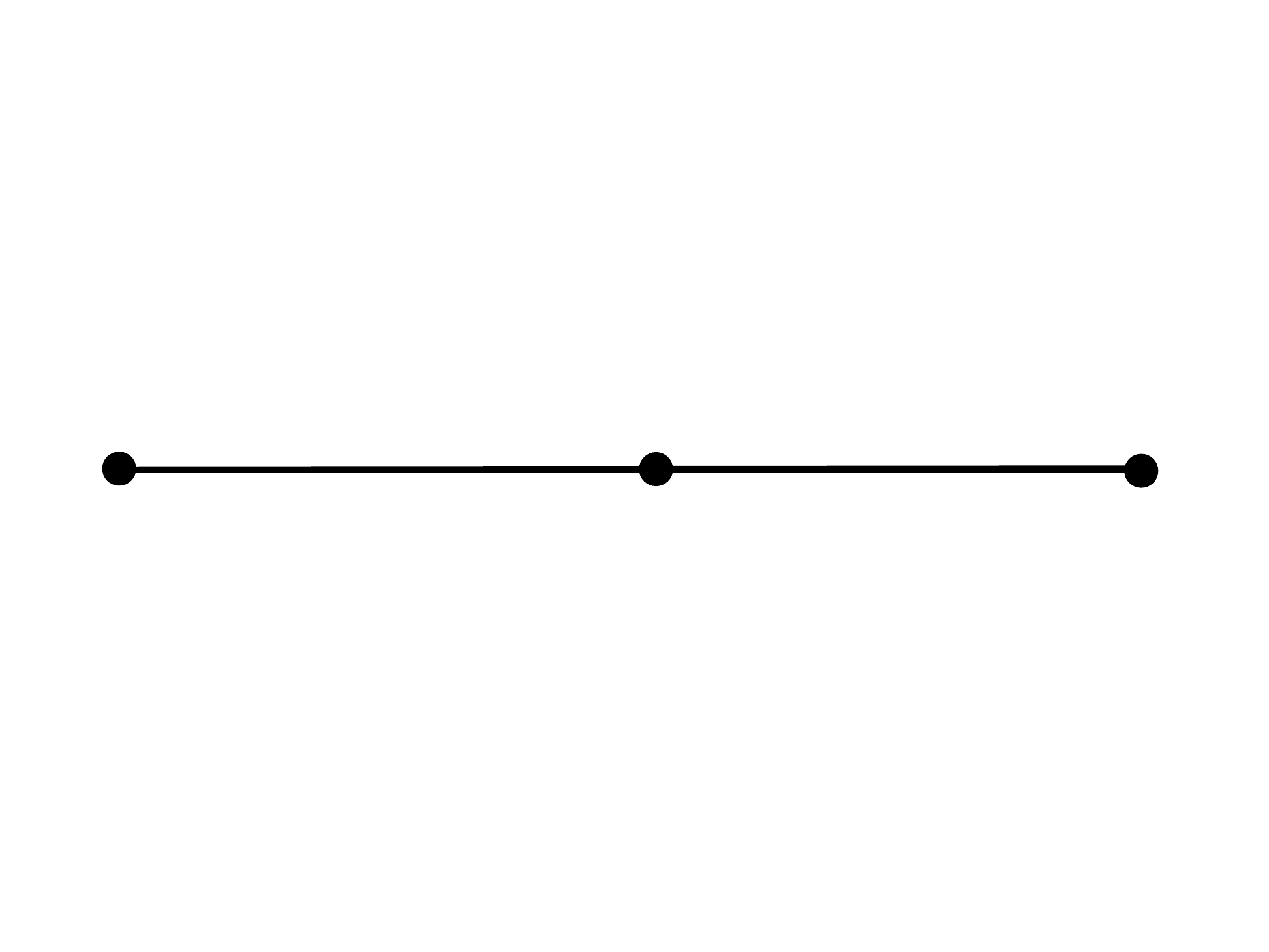}\hspace{1cm}\includegraphics[scale=.15]{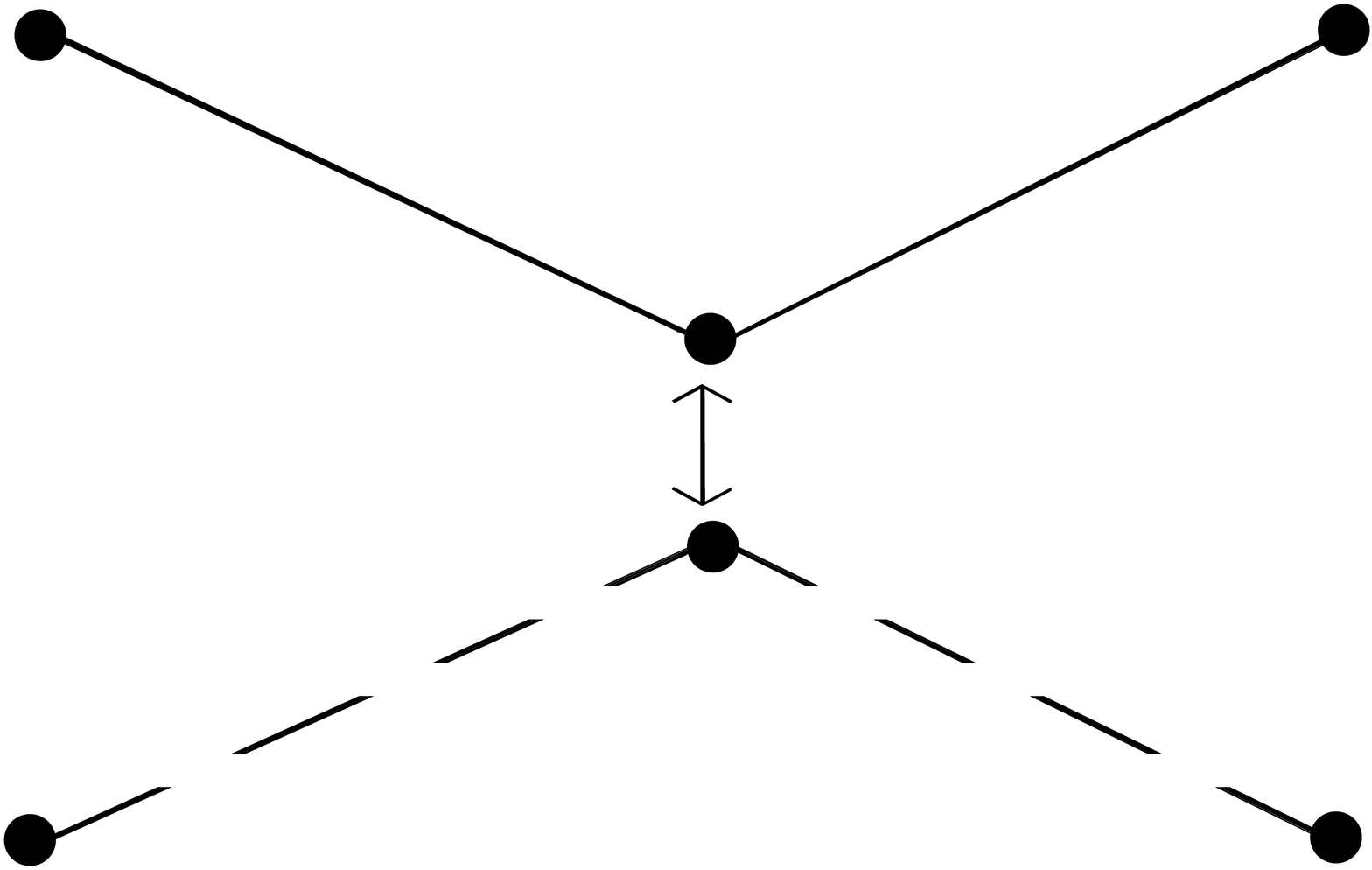}\\
$F_0$\hspace{4cm}$F_0\times\{0,1\}$\\ \vspace{1.5cm}
\includegraphics[scale=.15]{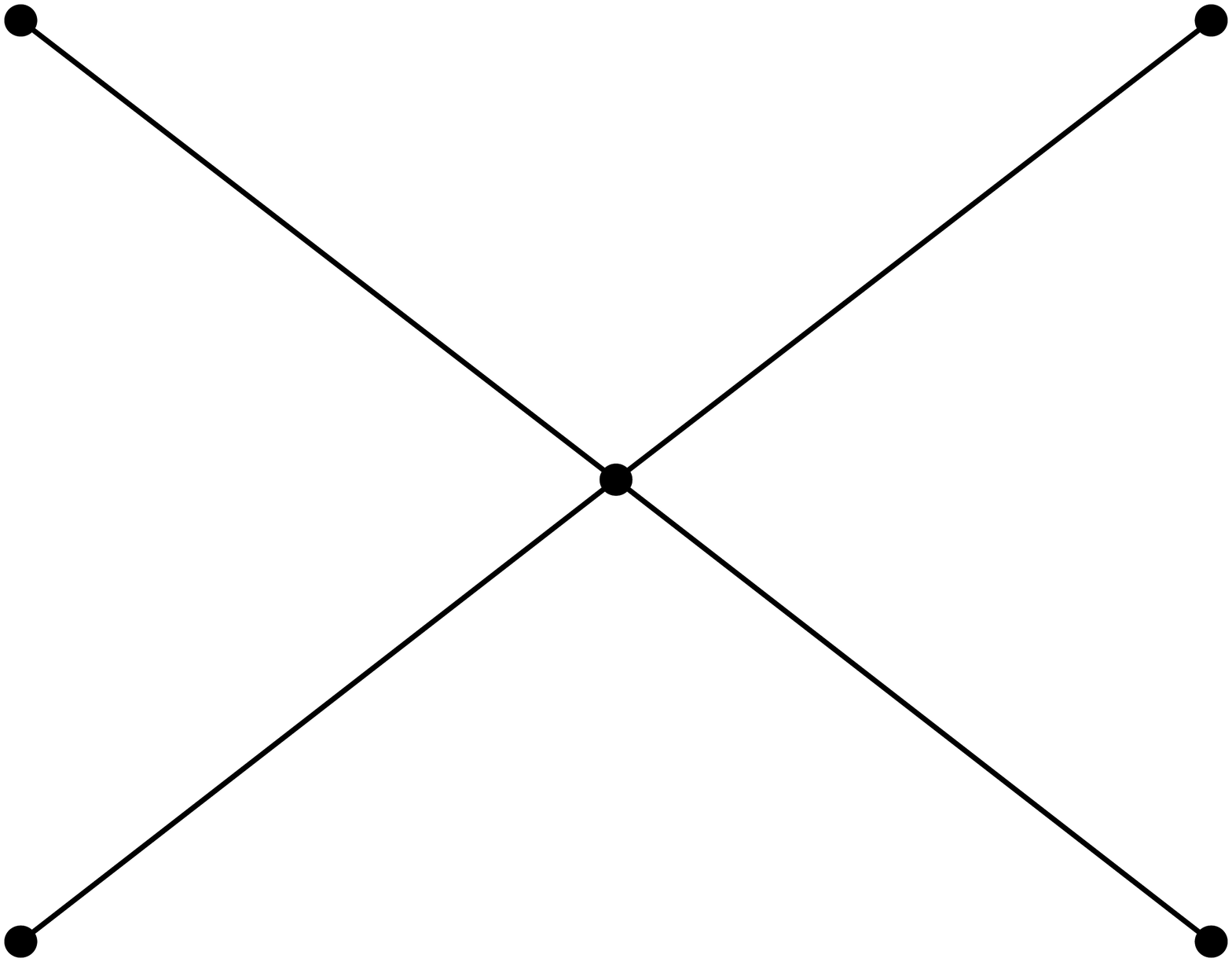}\\
$F_1$\\
\end{center}
\caption{Construction of $F_1$ from $F_0$ with $j_1=2$. }
\label{j2}
\end{figure}

We describe a simple case, where $j_n = 2$, $n\geq1$.  To obtain $F_1$, bisect $F_0$ with a node.  Then make a copy of this graph.  Identify the new nodes ``glueing" the two graphs together, represented by the arrow in Figure \ref{j2}.  This glueing process is the identification process described in \cite{BarlowEvans2004}.  This yields the graph $F_1$, an X-shape with five nodes as seen in Figure \ref{j2}. 

\begin{figure}[htbp]
\begin{center}
\includegraphics[scale=.15]{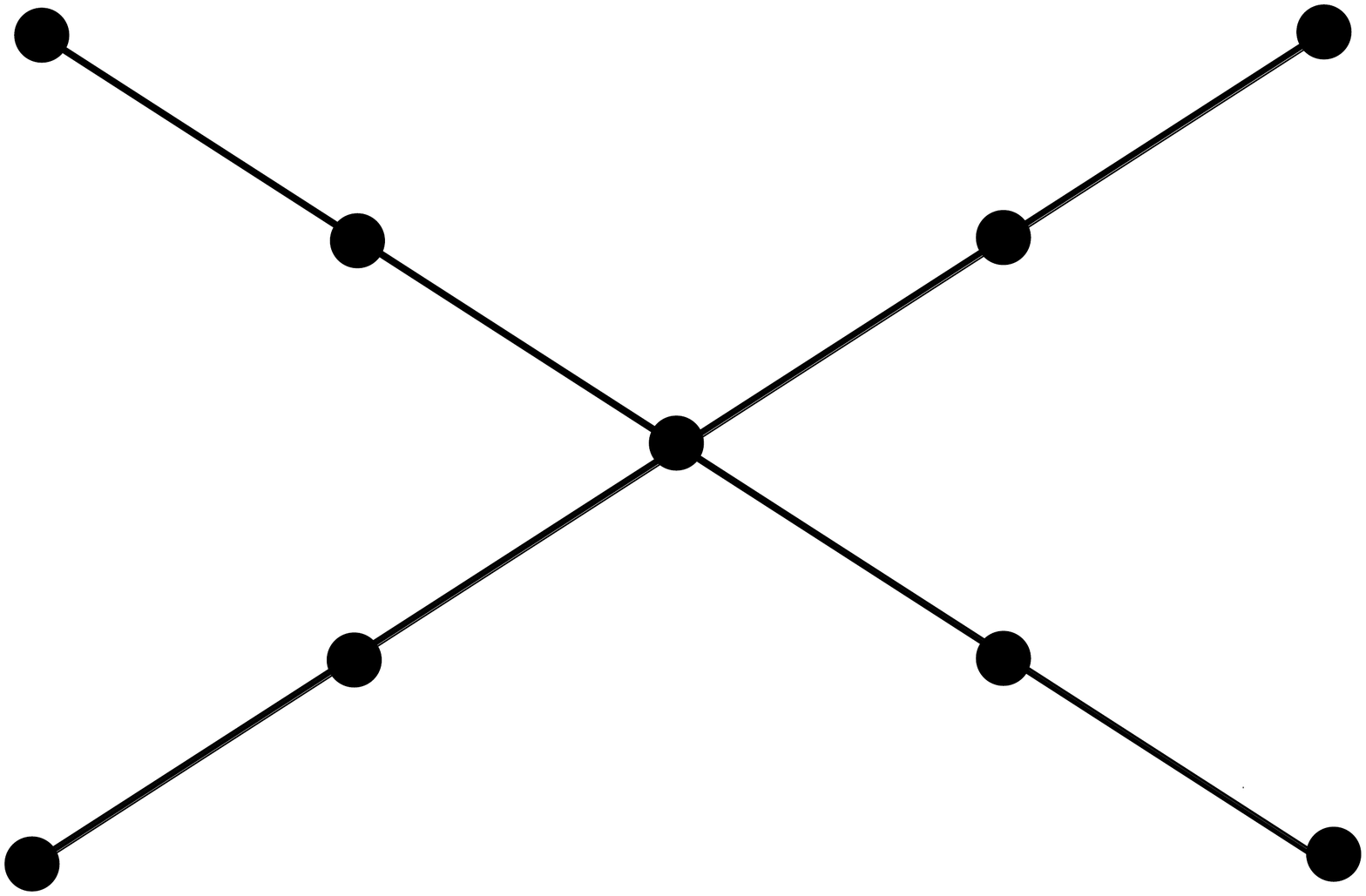}\hspace{2cm}\includegraphics[scale=.15]{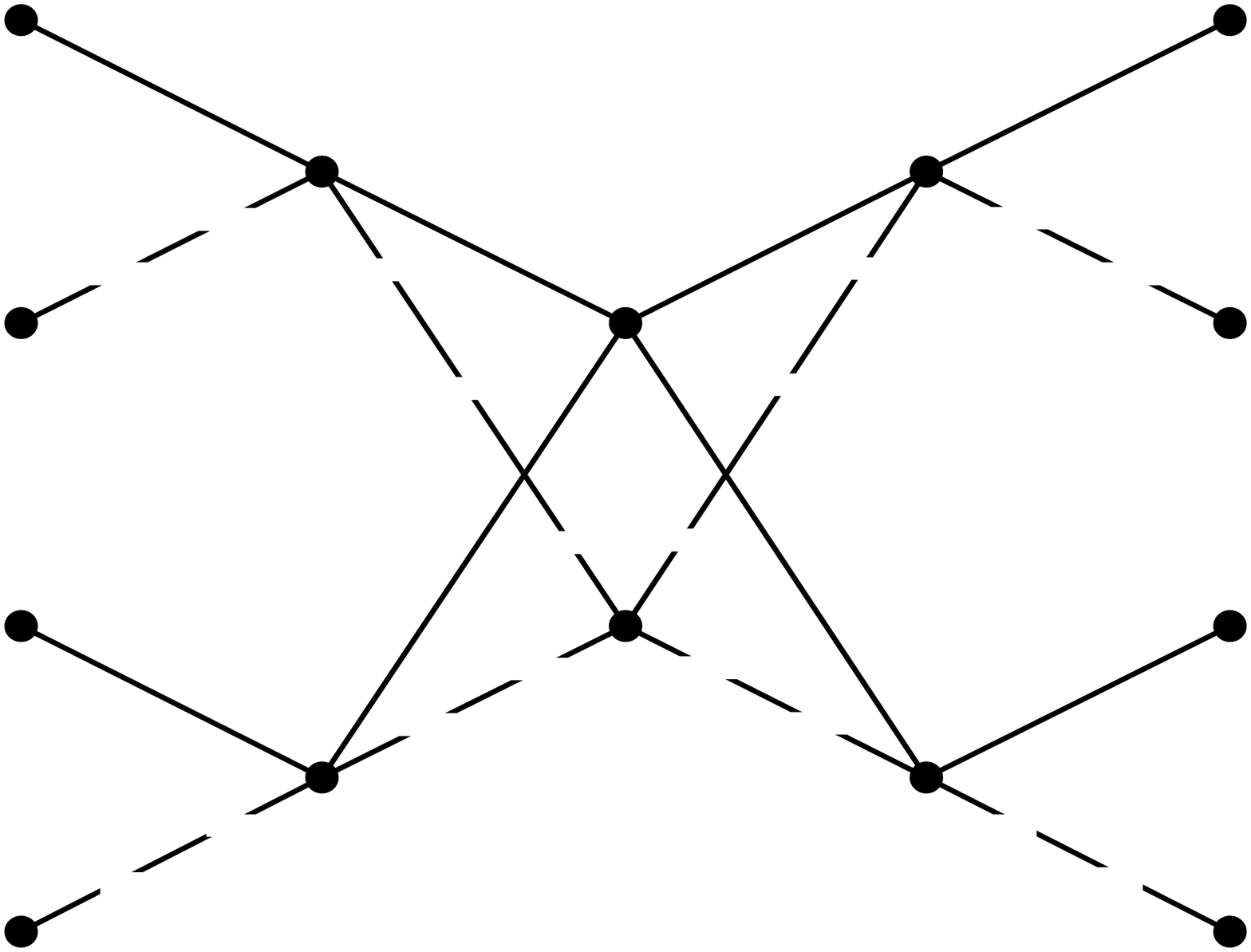}\\
$F_1$\hspace{5cm}$F_2$
\end{center}
\caption{Construction of $F_2$ from $F_1$ where $j_2=2$.  The dashed lines represent the second copy of $F_1$ with the added nodes.}
\label{j20}
\end{figure}

This procedure is repeated to obtain $F_2$ from $F_1$.  Nodes bisect each cell of $F_1$ as seen in Figure \ref{j20}.  A duplicate copy of $F_1$ is created and the two graphs are ``glued" together at the newly added nodes.  This is shown in Figure \ref{j20} where the solid line represents $F_1$ and the dashed lines represent the copy of $F_1$.  The $j_n=2$ Laakso space is the projective limit of the sequence of graphs $\{F_n\}_{n=0}^\infty$ all produced in this manner.

\begin{figure}[htbp]
\begin{center}
\includegraphics[scale=.2]{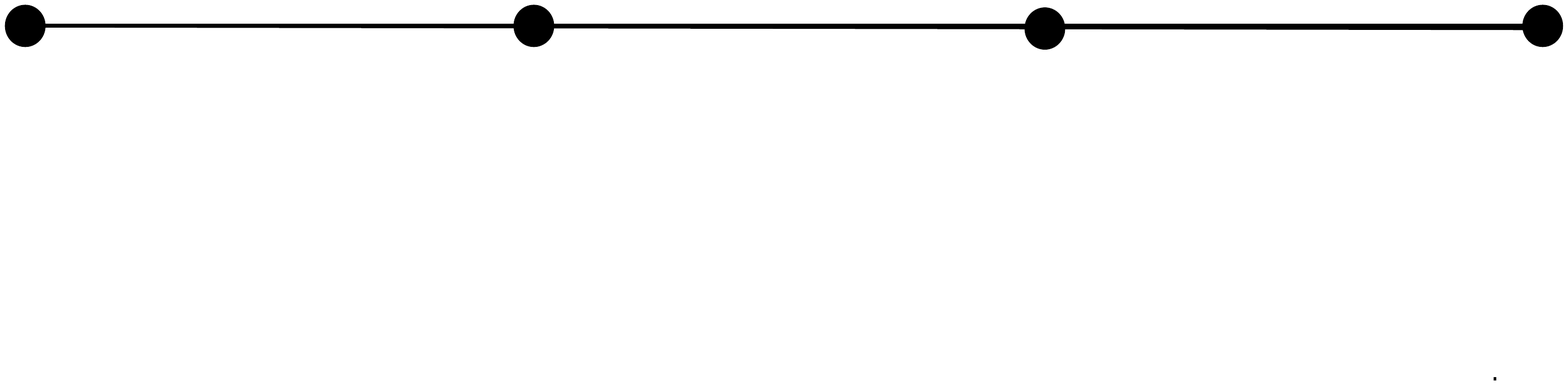}\includegraphics[scale=.2]{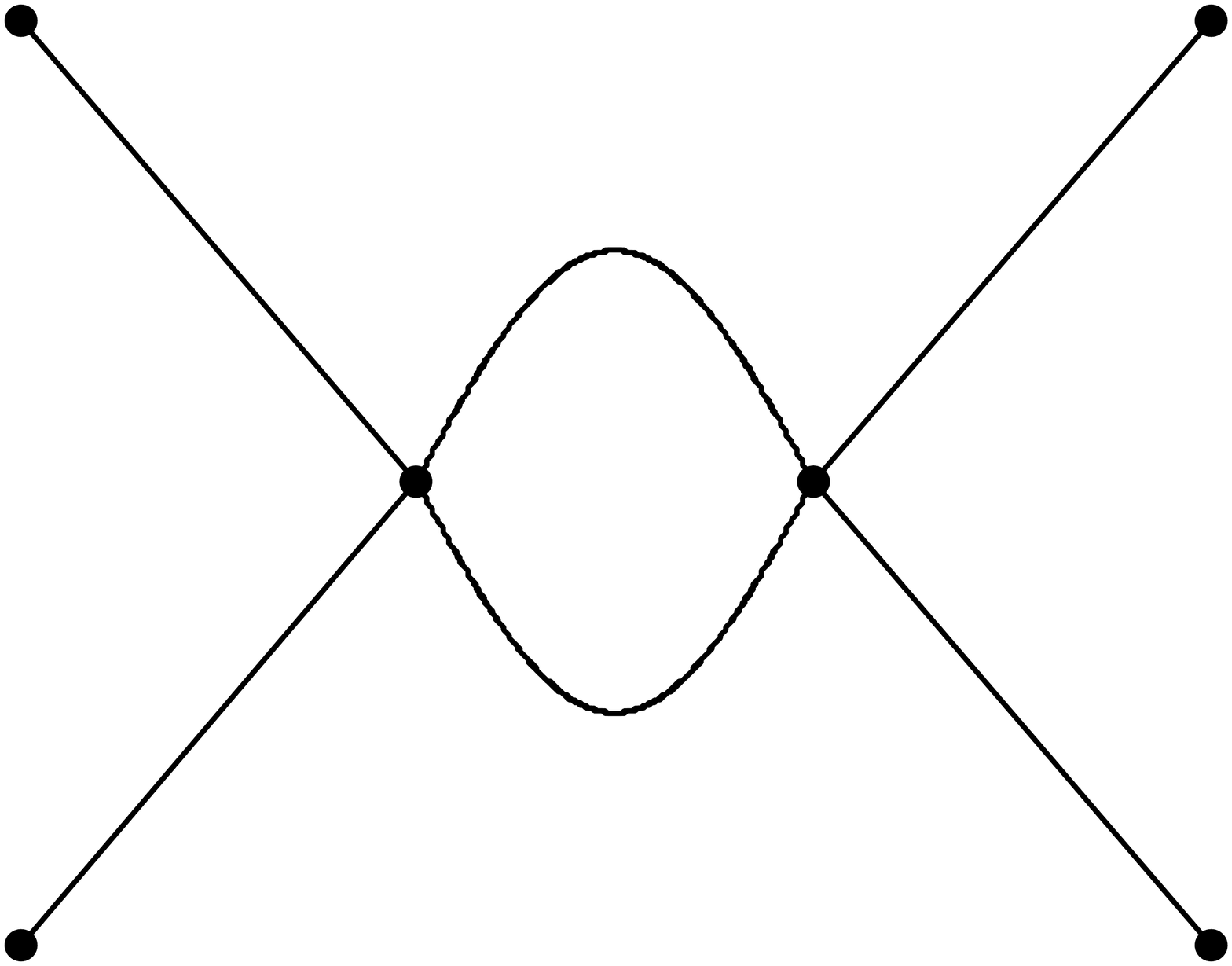}\\
$F_0$\hspace{5cm}$F_1$
\end{center}
\caption{Construction of $F_1$ from $F_0$ with $j_1=3$.}
\label{pic:j3}
\end{figure}

As another simple example consider $j_n = 3$,  for all $n\in \mathbb{N}$.  Again starting with the unit interval, $F_0$, $F_1$ is constructed by splitting $F_0$ into three subintervals and placing a node between each interval as shown in Figure \ref{pic:j3}.  This graph is duplicated and is glued to the original graph at the newly added nodes.  The two nodes in the middle of the figure are connected by the middle interval and its copy, thus creating a loop shape that is not seen in the case where $j=2$.  The outer thirds of the figure create a ``V" shape, also seen in the $j=2$ construction.  These shapes, loop and ``V", will be two of those considered in Section \ref{sec:delta}.

\begin{figure}[t]
\begin{center}
\includegraphics[scale=.2]{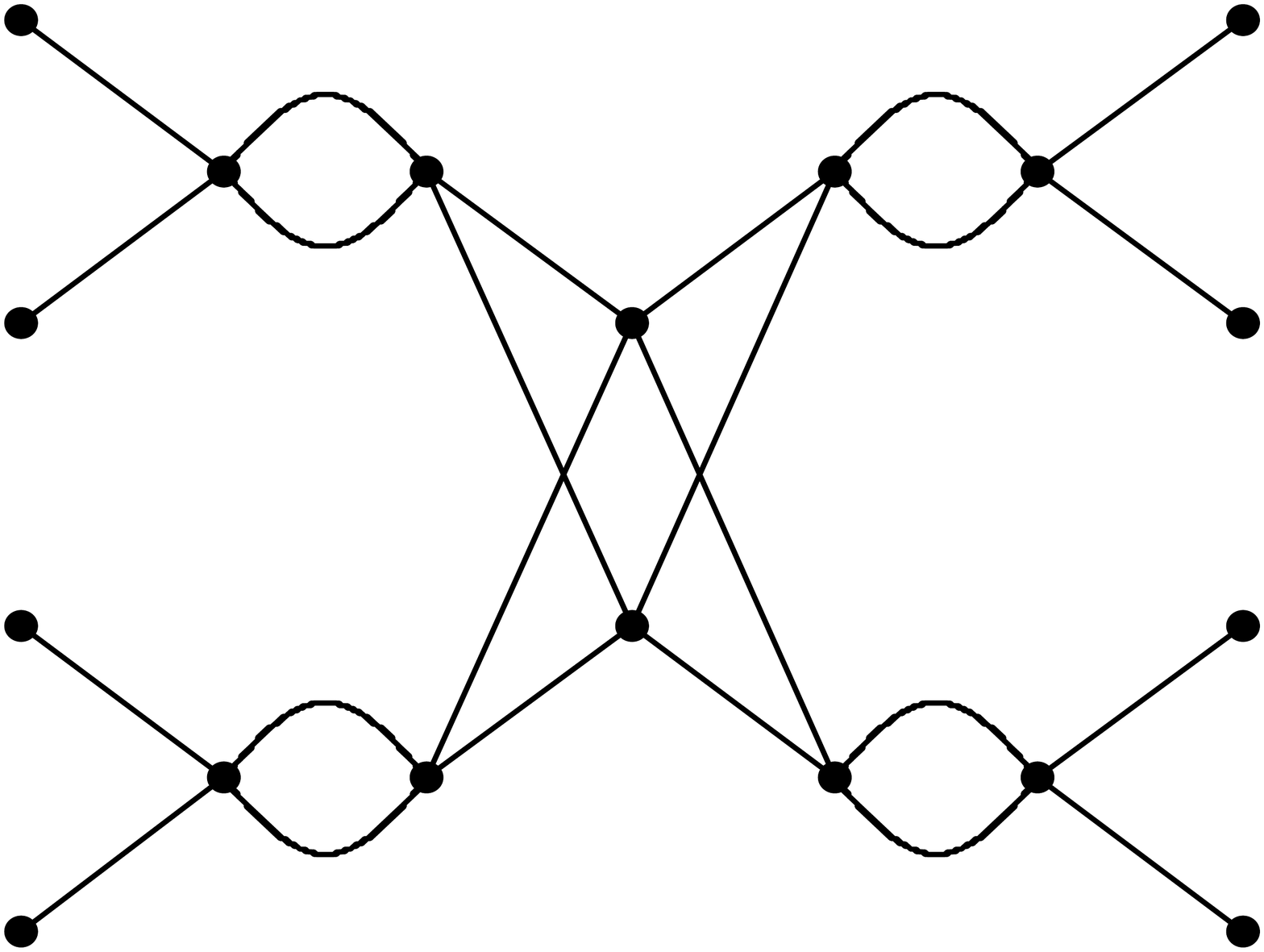}\includegraphics[scale=.2]{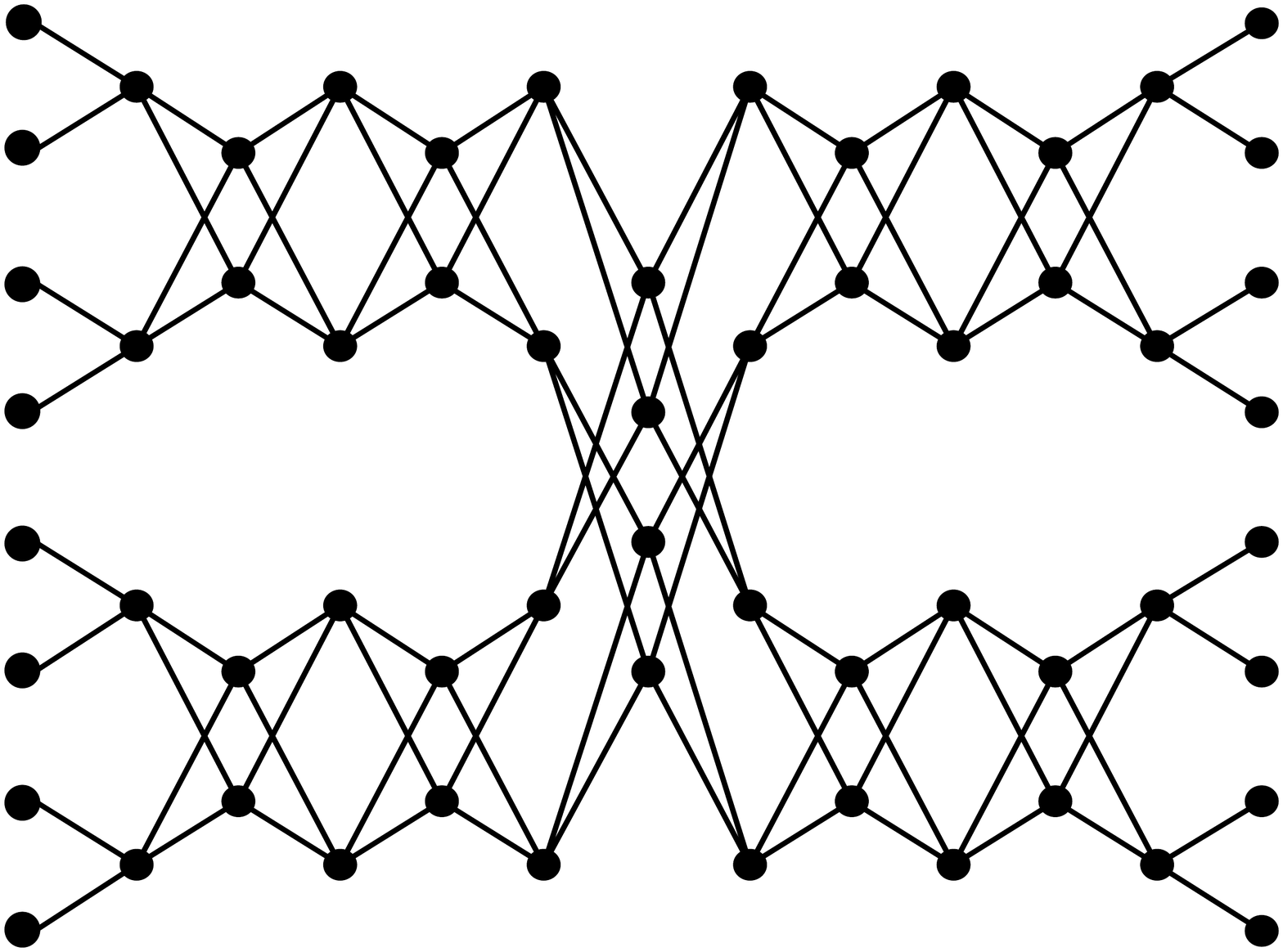}\\
$j_2$=[2,3] \hspace{5cm} $j_3=[2,3,2]$\\
\end{center}
\caption{Constructions of the Laakso space for $\{j_n\}_{n=1}^\infty=\{2,3,2,3,...\}$}\label{pic:2323}
\end{figure}

In this paper, we deal with the general case where $j_n$ may vary at each approximation level $n$.  The sequence $\{j_n\}_{n=1}^\infty$ may be a constant integer, as seen in the previous examples.  Or the sequence may alternate regularly between two integers.  Figure \ref{pic:2323} shows the construction of $F_2$ and $F_3$ when $\{j_n\}_{n=1}^\infty=\{2,3,2,3,...\}$.  The sequence $\{j_n\}_{n=1}^\infty$ could even be a completely random sequence of integers.  In any case, it is $\{j_n\}_{n=1}^\infty$ which defines the Laakso space and from which the properties are derived.

\subsection{Cell Structure of a Laakso Space}
Recall that as an inverse limit system the pair $(L, \{F_n\})$ come with continuous projection $\Phi_n:L \rightarrow F_n$.

\begin{defn}
The cell structure in a Laakso space, $L$, is determined by the pre-images under the map $\Phi_n$ of the cells in the graph $F_n$ which approximates $L$
\end{defn}
  
  Given a space as defined by Laakso in \cite{Laakso}, and the construction of the space as given in \cite{eigen}, and the level of approximation, $n$, the cell structure has specific properties, including number of cells,  $N_n$, and the interval length, $I_n^{-1}$.  Both the number of cells and the interval length are dependent on the choice of $j_i$ for all $i\leq n$. 

\begin{prop}\label{prop:cells}
Each cell in $F_n$ has metric diameter
\begin{equation}
I_n=\displaystyle\prod_{i=1}^nj_i
\end{equation}
where $I_0=1$. In addition the number of cells is
\begin{equation}
N_n=2^n\displaystyle\prod_{i=1}^nj_i.
\end{equation}
\end{prop}

\begin{proof}
In $F_0$ there is a single cell, the unit interval, with metric diameter equal to 1. At each step in the construction the diameter of each cell in $F_n$ is $j_n^{-1}$ times that of a cell in $F_{n-1}$. By induction the diameter of the cells in $F_n$ is $I_n^{-1}$. 

There is a single cell in $F_0$. At each step of the construction there are $2 \times j_n$ cells in $F_n$ for every cell in $F_{n-1}$. Thus the number of cells in $F_n$ is $2^{n}I_n$.
\end{proof}

\subsection{Hausdorff dimension of the Laakso Space}\label{hdim}  

In order to discuss the Hausdorff dimension of Laakso spaces we fix our choice of metric and measure. We use the path length metric. The measure used is the probability measure that gives equal mass to all cells of a given depth. Implicitly in the given construction, we have restricted the Hausdorff dimension to $1\leq Q \leq 2$.  In \cite{eigen} the Hausdorff dimension $Q$ of the Laakso Space associated with a constant $j_n$ at every level $n$ is shown to be
\begin{equation}
Q=1+\frac{log(2)}{log(j)}\nonumber
\end{equation}
Here we give the Hausdorff dimension of a Laakso space associated with a general sequence $\{j_n\}_{n=1}^\infty$. The measure used in calculating the Hausdorff dimension is the projective limit of Lebesgue measure on $F_n$ scaled to have total mass one for all $n$.

\begin{lemma}
Given sequence $\{j_i\}_{i=1}^\infty$ the Hausdorff Dimension, $Q$,  of the corresponding Laakso space is given by
\begin{equation}
Q_{j_i}=\displaystyle\lim_{n\to\infty}\frac{log\left(2^n\displaystyle\prod_{i=1}^nj_i\right)}{log\left(\displaystyle\prod_{i=1}^nj_i \right)}=\displaystyle\lim_{n\to\infty}\frac{log(2^nI_n)}{log(I_n)}=\lim_{n\to\infty} 1+\frac{log(2^n)}{log(I_n)},
\end{equation}
if the limit exists.
\end{lemma}

\begin{proof}
Laakso spaces are lacunary self-similar sets as defined in \cite{Igudesman2003} where the contraction ratios at any $n$ are equal. The number of identifications for each cell at the $i$'th iteration is $j_i$ and the formula that Igudesman gives in \cite{Igudesman2003} can be given in terms of $n$ and $j_i$.  This formula uses the number of cells, $N_n$ and the cell diameter, which is simply $I_n^{-1}$.  In the geodesic metric, the cell length is also the diameter of the cell.  The resulting formula is given above.
\end{proof}

While there are many sequences $\{j_n\}$ for which $Q$ will not exist it is more relevant to our interests that for every $Q$ there exist sequence $\{j_n\}$ yielding a Hausdorff dimension of $Q$. Different sequences $\{j_n\}$ can yield the same dimension, as shown in Table \ref{hdim:}.  These values agree with the dimensions given implicitly in \cite{Laakso}.

\begin{table}[htbp]\center
$\begin{array}{|c|c|}
\hline
$$j_i$$&$$Q_{j_i}$$\\ \hline
2&2\\  \hline
3&\frac{log(6)}{log(4)}\\ \hline
[2,3,2,3,...]&\frac{log(24)}{log(6)}\\ \hline
[3,2,3,2,...]&\frac{log(24)}{log(6)}\\ \hline
\end{array}$
\caption{Hausdorff Dimension for Laakso Space associated with given sequence of $\{j_i\}_{i=1}^\infty$}
\label{hdim:}
\end{table}

\subsection{Laplacian}\label{sec:Laplacian}
In \cite{eigen,ben} Laakso spaces are described as projective limits of quantum graphs and it is shown how to extend a compatible family of self-adjoint operators on the approximating quantum graphs to a self-adjoint operator on the limit space, i.e. the Laakso space. It was also shown how to use the spectrum with multiplicities of the operators on each quantum graph to determine the spectrum with multiplicities of the operator on the limit space. A quantum graph is a metric graph with a Hamiltonian operator,  as described in \cite{Kuchment2004,Kuchment2005}, the simplest of which would be the Laplacian operator, i.e.: a Hamiltonian without a potential.  

On each metric graph, $F_n$, consider the space of functions defined on the collection of edges each treated as a line segment. Define an operator on these function by $\Delta_n = -\frac{d^{2}}{dx^{2}}$. To make this a self-adjoint operator we need to also specify a suitable domain. A function is in $Dom(\Delta_n)$ if it is continuous everywhere, continuously twice differentiable on each line segment, and has Kirchhoff matching conditions at the nodes.  Kirchhoff conditions require a  function's with directional first derivatives summing to zero at nodes.

\begin{defn}
A Laakso space, $L$, is a projective limit of the $F_n$. Also there exist projection maps $\Phi_n:L \rightarrow F_n$ for all $n$. Thus any function on $F_n$ can bepulled back to a function on $L$ by writing $f \circ \Phi_n = \tilde{f}:L \rightarrow \mathbb{R}$. The pulling back is under the projections $\Phi_n$.
\end{defn}

By Theorem 7.1 in \cite{ben} those functions in $Dom(\Delta)$ that are pull backs are dense. so a complete set of eigenfunctions can be taken from this set. A consequence of this is that we can numerically approximate the spectrum of the Laplacian on the Laakso spaces working on some $F_n$.  Computations of these approximations are described in Section \ref{ssec:matlab} along with calculations described in Section \ref{counts}.  Tables \ref{table23} and \ref{23table} show calculated values of the spectrum of the Laplacian for Laakso space.

\begin{table}[tbh]\center\small
$\begin{array}{|cc|cc|cc|cc|cc||c|}\hline
n = 3 & & n = 4 & & n = 5 & & n = 6 & & n = 7 &&  Expected    \\
\lambda & m & \lambda & m  & \lambda & m & \lambda & m & \lambda & m  &   \\ \hline
0 & 1 & 0 & 1 & 0 & 1 & 0 & 1& 0 & 1&0    \\
9.87 & 3 & 9.87 & 3 & 9.87 & 3 & 9.87 & 3 & 9.87 & 3 & \pi^2 = 9.87 \\
39.58 & 1 & 39.38 & 1 & 39.45 & 1 & 39.48 & 1 & 39.48 & 1& (2\pi)^2 = 39.48\\
84.35 & 8 & 88.32 & 8 & 88.70 & 8 & 88.81 & 8 &88.82 & 8 & (3\pi)^2 = 88.83\\
144 & 1 & 156.32 & 1 & 157.51 & 1 & 157.87 & 1 & 157.90& 1 & (4\pi)^2 = 157.91\\
 213.46& 3 & 242.85 & 3 & 245.76 & 3 & 246.63 & 3 & 246.71 & 3 & (5\pi)^2 = 246.74\\
 288 & 26 & 347.26 & 26 & 353.28 & 26 & 355.08 & 26 & 355.25 & 26 & (6\pi)^2 = 355.31\\
 & & 468.78 & 3 & 479.86 & 3 & 483.19 & 3 & 483.51 & 3 & (7\pi)^2 = 483.61\\
 & & 606.41 & 1 & 625.27 & 1 & 630.94& 1 & 631.48  & 1 & (8\pi)^2 = 631.65\\
 & & 759.18 & 8  & 789.22 & 8 & 798.30 & 8 & 799.15  & 8 &  (9\pi)^2 = 799.44\\
 & & 925.89 & 1  & 971.40 & 1 & 985.22 & 1 & 986.53 & 1 &  (10\pi)^2 = 986.96\\
 & & 1105.3 & 3  & 1171.5 & 3 & 1191.7 & 3 & 1193.6 & 3 & (11\pi)^2 = 1194.2\\
 & & 1296 & 38  & 1389.0 & 38 & 1417.6 & 38 & 1420.3 & 38 & (12\pi)^2 = 1421.2\\
 & & 1496.6 & 3  & 1623.7 & 3 & 1663.0 & 3 & 1666.7 & 3 &  (13\pi)^2 = 1668.0\\
 & & 1705.5 & 1  & 1875.0 & 1 & 1927.8 & 1 & 1932.8 & 1 &  (14\pi)^2 = 1934.4\\
 & & 1921.1 & 8  & 2142.5 & 8 & 2211.9 & 8 & 2218.5 & 8 &  (15\pi)^2 = 2220.7\\
 & & 2141.9 & 1  & 2425.7 & 1 & 2515.2 & 1 & 2523.8 & 1 &  (16\pi)^2 = 2526.6\\
 & & 2366.1 & 3  & 2723.9 & 3 & 2837.8 & 3 & 2848.7 & 3 &  (17\pi)^2 = 2852.3\\
 & & 2592 & 86  & 3036.7 & 86 & 3179.5 & 86 & 3193.2 & 86 &  (18\pi)^2 = 3197.8\\
 & & 2817.9 & 3 & 3363.5 & 3 & 3540.3 & 3 & 3557.3 & 3 & (19\pi)^2 = 3562.9\\
 \hline
 \end{array}$
 \caption{Calculated Values of the first 20 Eigenvalues for $\{j_n\}_{n=1}^\infty=\{2,3,2,3,...\}$ with multiplicity, $m$, the iteration value, $n$, and the expected value, $\lambda$.  As $n$ increases, the observed eigenvalues converge to the expected result.}
 \label{table23}
 \end{table}

\section{The Spectrum of $\Delta$}\label{sec:delta}
Theorem \ref{thm:spect} gives the spectrum and associated multiplicities of the Laplacian operator by considering $\Delta_n$ on $F_n$ and on any Laakso space.  We devote the rest of the section to proving the theorem.  Following the analytic arguments are details of computational experiments carried out before the analytic results were available.  We use an iterative, computer-assisted process to find the bottom end of the spectrum on a number of specific Laakso spaces.  In all cases the computed results and analytic results agree within the precision of the computations.    

\begin{theorem}\label{thm:spect}
Given any Laakso space, $L$, with associated sequence $\{j_i\}_{i=1}^n$, the spectrum of $\Delta$ on $Dom(\Delta)$ is 
\begin{eqnarray}\small
\bigcup_{k=0}^{\infty} \{\pi^2 k^2\}\cup\bigcup_{n=1}^\infty \bigcup_{k=0}^{\infty} \{(k+1/2)^2\pi^2 I_n^2\}\cup\bigcup_{n=1}^\infty \bigcup_{k=1}^\infty \{k^2\pi^2I_n^2\}\notag\\ \cup\bigcup_{n=2}^\infty \bigcup_{k=1}^\infty \{k^2\pi^2I_n^2\} \cup \bigcup_{n=2}^\infty \bigcup_{k=1}^{\infty} \left\{\frac{k^2\pi^2 I_n^2}{4}\right\}\notag\\
\end{eqnarray}
with associated multiplicities:
\begin{equation}
1,\hspace{.5cm} 2^n,\hspace{.5cm} 2^{n-1}(j_n-2)I_{n-1},\hspace{.5cm} 2^{n-1}(I_{n-1}-1),\hspace{.5cm} 2^{n-2}(I_{n-1}-1)
\end{equation}\label{eqn:multi}
respectively.
\end{theorem}

This does correct a typographical error in the similar statement given in \cite{eigen}.

\begin{table}[t]\center\small
$\begin{array}{|cc|cc|cc|cc||c|}
\hline
$A$ & & $B$ & & $C$ & & $D$& &   Expected    \\
\lambda & m & \lambda & m  & \lambda & m & \lambda & m & \\ \hline
0 & 1 & 0 & 1 & 0 & 1 & 0 & 1&0    \\
9.87 & 3 & 9.87 &1 & 9.87 & 1 & 9.87 &1  &   \pi^2 \\
& & 22.2 & 2 & 22.2 & 2 &  &  &   (1.5\pi)^2\\
39.48 & 1 & 39.48 & 1 & 39.48 & 1 &39.48 &3  &  (2\pi)^2\\
88.82 & 8 & 88.82 & 8 & 88.82 & 2 &88.82 &1 &  (3\pi)^2\\
157.91& 1 & 157.91 & 1 & 157.91 & 1 &157.91  &3  &  (4\pi)^2\\
& & 199.85.3 & 2 & 199.85 & 2 & &  &   (4.5\pi)^2\\
246.71 & 3 & 246.71 & 1 & 246.71 & 1 & 246.71 &1  & (5\pi)^2\\
355.25 & 26 & 355.25 & 8 & 355.25 & 8 &355.25 & 10 &  (6\pi)^2\\
483.51  & 3& 483.51 & 1 & 483.51 & 1 & 483.51 & 1 & (7\pi)^2\\
& & 555.1 & 2 & 555.1 & 2 &  &  &  (7.5\pi)^2\\
 631.48& 1& 631.48 & 1 & 631.48 & 1 &631.48 &3  &    (8\pi)^2\\
 799.15& 8& 799.15 & 36  & 799.15 & 2 & 799.15 & 1 &    (9\pi)^2\\
 986.53& 1& 986.53 & 1  & 986.53 & 1 & 986.53& 3 &   (10\pi)^2\\
 &  & 1087.9& 2 & 1087.9& 2 & &  &   (10.5\pi)^2\\ 
 1193.6&3 & 1193.3 & 1  & 1193.3& 1 &1193.3  &1  &   (11\pi)^2\\
 1420.3&38 & 1420.3 & 8  & 1420 & 2 &1420 &20 &  (12\pi)^2\\
 1666.7& 3& 1666.7& 4  & 1666.7 & 1 &1666.7 &1 &   (13\pi)^2\\
 &  & 1798.1& 2 & 1798.1 & 2 & &  &  (13.5\pi)^2\\ 
 1932.8&1 & 1932.8& 1  & 1932.8& 1 &1932.8 &3  &   (14\pi)^2\\
 2218.5&8 & & & &  & 2218.5 & 1 & (15\pi)^2\\
 2523.8&1 & & & &  & 2523.8 & 3 &    (16\pi)^2\\
 2848.7&3 & & & & &2848.7 &1 &  (17\pi)^2\\
 3193.2&86 & & & &  &3139.2 &10  &  (18\pi)^2\\
3557.3 & 3& & & &&3557.3 &1& (19\pi)^2\\
\hline
 \end{array}$
 \caption{Calculated Values of the first 20 eigenvalues for given sequences of $j_i$'s with multiplicity, $m$ and the expected value, $\lambda$.}
 A=\{2,3,2,3,...\}\hspace{0.5cm} B=\{3,2,3,2,...\} \hspace{0.5cm} C=\{3,4,3,4,...\}\hspace{0.5cm} D=\{4,3,4,3,...\}
 \label{23table}
 \end{table}
 
\subsection{Counts of Eigenvalues and Multiplicities}\label{counts}
\begin{figure}[tbp]
\begin{center}
\includegraphics[scale=.1]{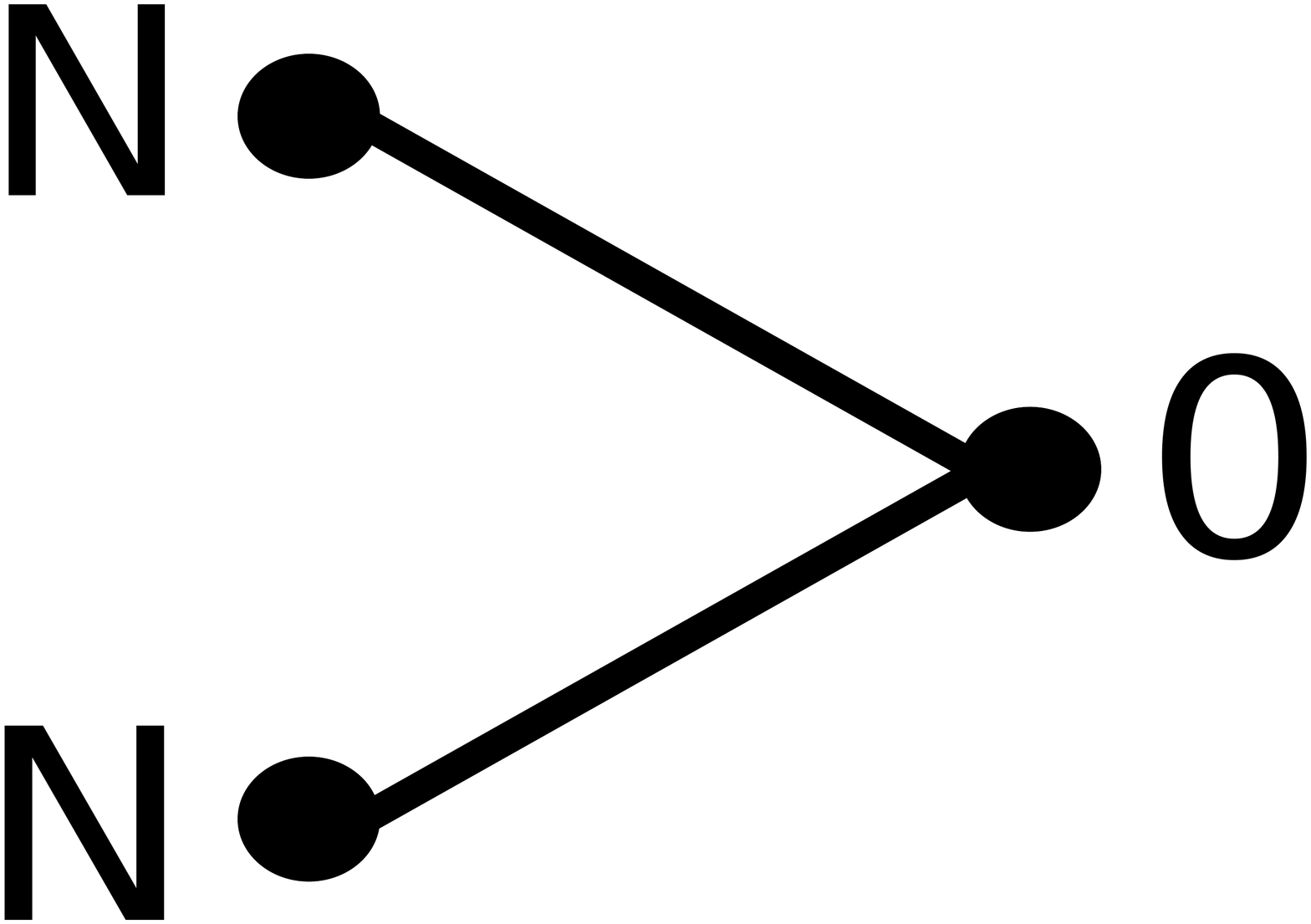}\hspace{1cm}\includegraphics[scale=.1]{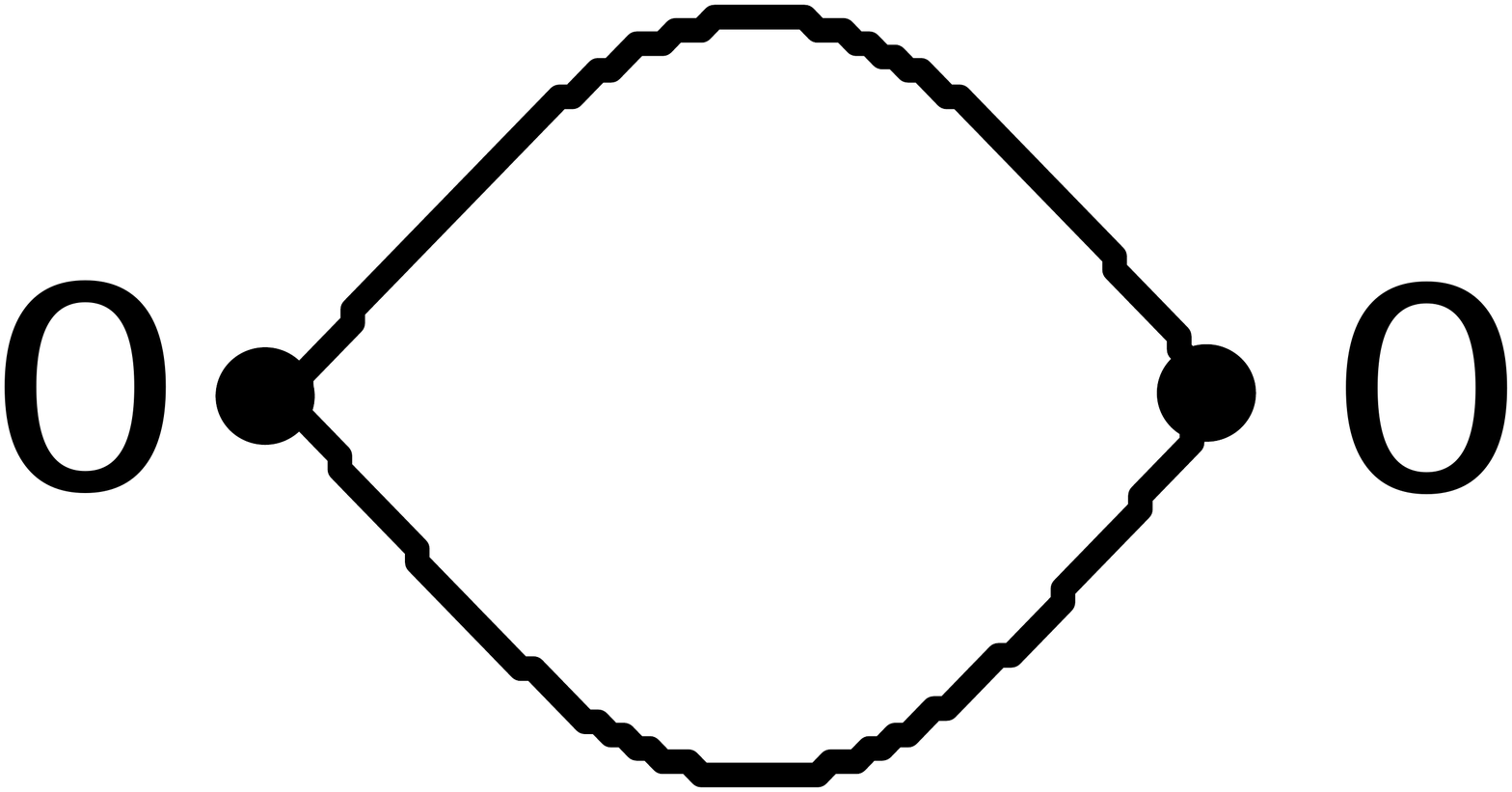}\hspace{1cm}\includegraphics[scale=.1]{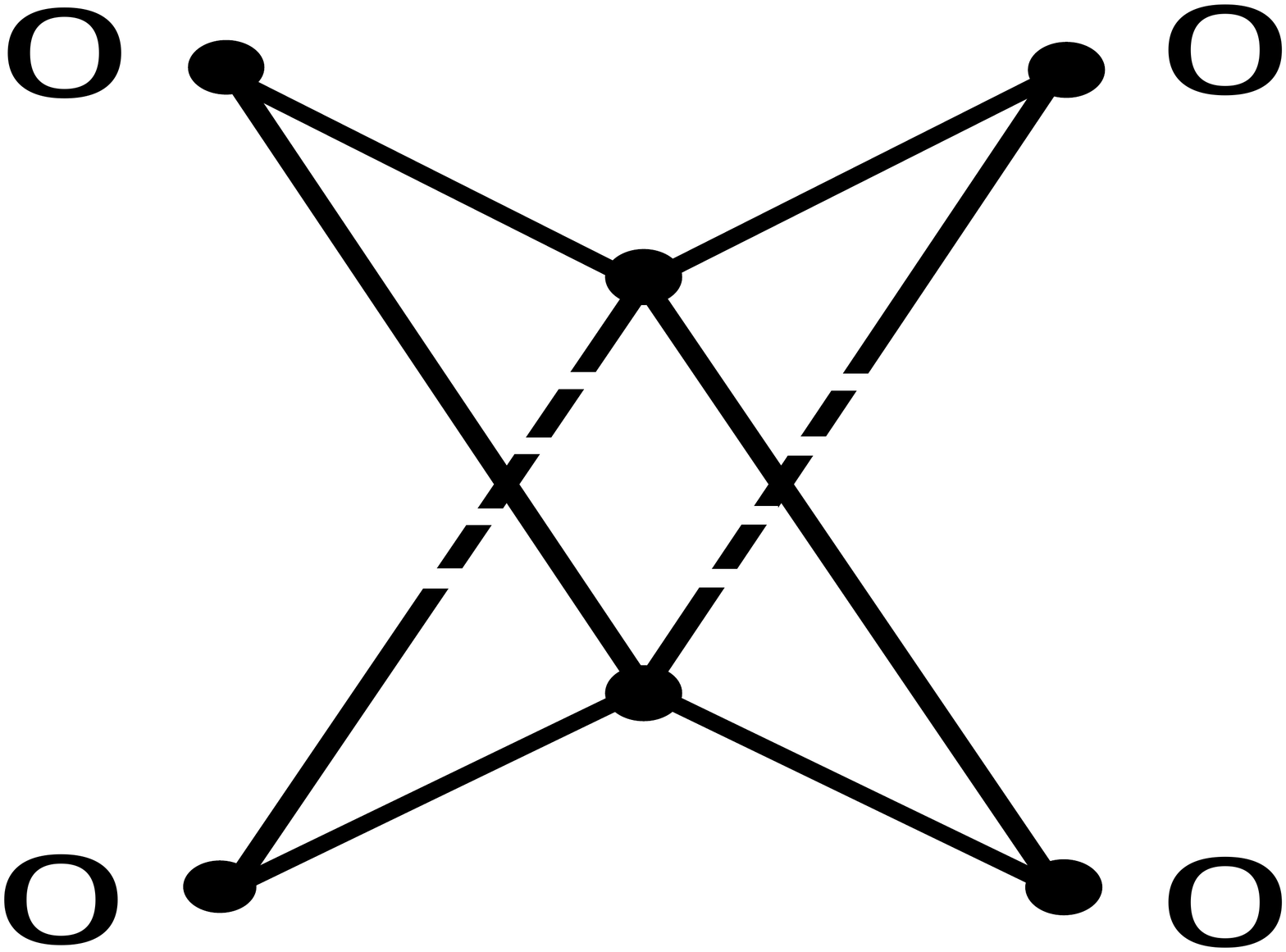}\\
$``V"$\hspace{3cm}$``Loop"$\hspace{3cm}$``Cross"$
\end{center}
\caption{Constructions of v's, loops, and crosses along with associated boundary conditions}
\label{pic:loop}
\end{figure}

In order to determine the spectrum of the Laplacian on the Laakso space, the approximating quantum graph is considered as a collection of simpler parts.  In \cite{eigen} it was determined that three distinct shapes with appropriate boundary conditions could be used to construct any quantum graph representations of Laakso spaces, save $F_0$, which is treated as a special case.  Definition \ref{shapedefn} defines these three shapes shown in Figure \ref{pic:loop} with their respective boundary conditions which are forced by the Kirchhoff matching conditions and the orthogonality requirements that assign an eigenfunction to a given representation level. These orthogonality conditions were discussed in detail in \cite{eigen}. In short they allow the counting arguments  to count an eigenfunction only once. 

\begin{defn}\label{shapedefn}
(a) A \emph{shape} is a connected quantum sub-graph, as shown in Figure \ref{pic:loop}. In that figure the ``D'' denotes a Dirichlet boundary condition at that node and an ``N''  the Neumann condition.

(b) A \emph{V} is the shape consisting of three nodes: two nodes in a column and the third node a second.  The two nodes in the first column are degree one and the node in the second column is degree two.  When a \emph{V} is in $F_n$, it shares its degree two node with another shape thus making it a degree four node, as seen in Figure \ref{pic:2323}.

(c) A \emph{loop} is the shape that consists of two nodes each of degree two.  The nodes are connected to each other by two cells, creating a loop. When a \emph{loop} is in $F_n$, both degree two nodes are shared as degree two nodes for another shape thus making them degree-four nodes, Figure \ref{pic:2323}. 

(d) A \emph{cross} is the shape consisting of six nodes four of degree two and two of degree four.  The degree two nodes each have a cell connecting the node to each of the degree four nodes .    Notice the cross is the only shape containing nodes of degree four in the subgraph.   When a \emph{cross} is in $F_n$ the degree two nodes are shared with degree two nodes of another shape thus making them degree four nodes, as in Figure \ref{pic:2323}. 
\end{defn}

Before determining the spectrum of $\Delta_n$  on the three shapes, we must first establish the following proposition which describes how these three shapes are involved in the construction of $F_n$ and $L$.

\begin{prop}\label{cellprop}
\begin{itemize}
	\item[(a)] Any node in any quantum graph approximating a Laakso space is either of degree one or degree four. 
	\item[(b)] For any degree one node in $F_n$, a \emph{V} is produced in $F_{n+1}$.  
	\item[(c)] For any degree four node in $F_n$, a \emph{cross} is produced in the construction of $F_{n+1}$.  
	\item[(d)] Any cell in $F_n$ produces $j_{n+1}-2$ \emph{loops} in $F_{n+1}$ between the \emph{V}'s or \emph{crosses} produced by the nodes in $F_n$.
	\item[(e)] For $n\geq0$ the number of nodes in $F_n$ is $N_n=2^{n-1}(I_n+3)$.
\end{itemize}
\end{prop}

\begin{proof}
\begin{itemize}
	\item[(a)] A degree one node in $F_{n-1}$ gives rise to two degree one nodes in $F_n$ as an immediate consequence of the construction. Similarly a degree four node gives rise to a single degree four node in $F_n$. The new nodes in $F_n$ that are not nodes in a copy of $F_{n-1}$ are the identification of two degree two nodes, hence of degree four. 
	\item[(b)] In the construction of $F_{n+1}$, the cell connected to a degree one node is split into $j_{n+1}$ intervals by adding $(j_{n+1}-1)$ nodes.  Then the graph is duplicated yielding two rows of cells connected between $j_{n+1}$ columns with two nodes in each, all of degree two except for the nodes at the end of the cell.  The original and duplicated cells are connected at the newly added nodes.  Thus the original node of degree one from $F_n$ remains degree one in $F_{n+1}$.  The graph around the original node and it's duplicate is a ``V.''  
	\item[(c)] In the construction of $F_{n+1}$, the four cells connected to the degree four node in $F_n$ will be split into $j_{n+1}$ intervals.  To construct $F_{n+1}$, $F_n$ is duplicated, new nodes inserted, and connected at the newly added nodes.  Thus the original degree-four node remains degree four and is duplicated, creating two degree-four nodes.  The graph around the original node and it's duplicate is a ``cross.''
	\item[(d)] Parts b and c account for two of the $j_{n+1}$ intervals.  The rest produce loops.  So, there are $j_{n+1}-2$ loops in $F_{n+1}$ for every cell in $F_n$.
	\item[(e)] We will induct on $n$.  The unit interval, $F_0$, has two nodes. Suppose that $F_{n-1}$ has $N_{n-1} = 2^{n-2}(I_{n-1}+3)$ nodes. Then $N_n = 2 \times N_{n-1} + 2^{n-1}(j_n-1)I_{n-1}$, the nodes from the two copies of $F_{n-1}$ plus the new nodes of which there are $j_n-1$ new nodes per cell in $F_{n-1}$ and there are $2^{n-1}I_{n-1}$ cells in $F_{n-1}$. This simplifies to the claimed formula.
\end{itemize}
\end{proof} 

We now generalize the results from \cite{eigen} in three lemmas that give the eigenvalues and multiplicites (counts) for each of the three shapes.

\begin{lemma}\label{vlemma}
For any $n\geq1$, the number of \emph{V}'s in $F_n$ is $2^n$.  The eigenvalues for this shape at this level are: 
\begin{equation}\label{vspect}
\{ [I_n (k+1/2) \pi]^2 : k=0,1,\ldots \}.
\end{equation}
\end{lemma} 

\begin{proof}
We prove the count by induction.  $F_1$ is constructed out of $F_0$, which is a single cell connecting two degree one nodes.  This implies by Proposition \ref{cellprop} that $F_1$ will have 2 \emph{V}'s.  Now assume that for some arbitrary $n\geq1$, the number of \emph{V}'s in $F_n$ is $2^n$.  From Definition \ref{shapedefn}, the \emph{V} is the only shape that has a degree one node.  Furthermore, it has two degree one nodes.  From Proposition \ref{cellprop} we know that each degree one node in $F_n$ produces a \emph{V} in $F_{n+1}$.  From \cite{eigen} the shapes defined in Definition \ref{shapedefn} are all the possible shapes in the graphs, so there cannot be any degree one nodes from any other shape.  So the number of \emph{V}'s in $F_{n+1}$ is twice the number of \emph{V}'s in $F_n$.  So $F_{n+1}$ has $2^{n+1}$ \emph{V}'s.

In order to get the spectrum of $\Delta_n$ restricted to a \emph{V} we look at the functions in this domain that are orthogonal to the functions expressible on an interval.  These functions have the property that the values on the top branch are the negative of the values on the lower.  We therefore need only consider the top branch, as it fully determines the behavior on the bottom branch.  This top branch is one interval, and has Neumann boundary conditions at one end and Dirichlet boundary conditions at the other.  The length of the cell is $I_n^{-1}$.  So we are looking for eigenfunctions on intervals of length $I_n^{-1}$ with zero derivative at one end and zero value at the other.  These come in the form $\cos(I_n(k+1/2)\pi x)$ where $k=0,1,\ldots$ and $x \in [0,I_n^{-1}]$.  The eigenvalues in \eqref{vspect} are now obtained in the usual way. 
\end{proof}

\begin{lemma}\label{looplemma}
For any $n\geq1$, the number of \emph{loops} in $F_n$ is 
\begin{equation}\label{loopcount}
2^{n-1}(j_n-2)(I_{n-1}).
\end{equation}
The eigenvalues for this shape at this level are
\begin{equation}\label{loopspect}
\{ [I_n k \pi]^2 : k=1,2,\ldots \}.
\end{equation}
\end{lemma} 

\begin{proof}
By Proposition \ref{cellprop} every cell in $F_{n-1}$ produces $j_{n}-2$ loops in $F_n$.  In order to know how many loops are in $F_n$, the number of cells in $F_{n-1}$ are counted and multiplied by $j_n-2$.  The number of cells in $F_n$ were already counted in Proposition \ref{prop:cells} and shown to be $2^n(I_n)$.  Substituting in $n-1$ for $n$ in this expression and multiplying by $j_n-2$ gives (\ref{loopcount}).  

In order to get the spectrum of $\Delta_n$ restricted to a loop we look at the functions in this domain that are orthogonal to the functions expressible on an interval.  Again, these functions have the property that the vales on the top branch are the negative of those on the lower.    As was the case with the \emph{V} above, the orthogonality condition imposed on the functions reduces the question to only considering the top interval of length $I_{n}^{-1}$ with Dirichlet boundary conditions.  The eigenfunctions that fit these conditions are $\sin(I_n k \pi x)$ with $k=1,2,\ldots$ and $x \in [0,I_n^{-1}]$.  These result in the set defined in (\ref{loopspect}).  

\end{proof}

\begin{lemma}\label{crosslemma}
For any $n\geq2$, the number of \emph{crosses} in $F_n$ is 
\begin{equation}\label{crosscount}
2^{n-2}(I_{n-1}-1).
\end{equation}
There are two sets of eigenvalues for this shape at this level.  They are 
\begin{equation}\label{crossspect1}
\left\{ \left[\frac{1}{2} (I_n k \pi)\right]^2 : k=1,2,\ldots \right\} 
\end{equation}
with multiplicity one and 
\begin{equation}\label{crossspect2}
\{ [I_n k \pi]^2 : k=1,2,\ldots \}
\end{equation}
with multiplicity two.
\end{lemma}

\begin{proof}
From Proposition \ref{cellprop} crosses in $F_n$ appear only where there were degree four nodes in $F_{n-1}$.  Therefore, to find the number of crosses in $F_n$, we will count the number of degree four nodes in $F_{n-1}$.  By Proposition \ref{cellprop} every node in a quantum graph approximating a Laakso space is either of degree one or degree four.  Therefore, subtracting the number of degree one nodes from the total number of nodes will give the number of degree four nodes.  From the same proposition, the total number of nodes is $2^{n-1}(I_n+3)$.  We have seen already that in $F_n$ degree one nodes only occur in \emph{V}'s and that for every v there are two degree one nodes.  From Lemma \ref{vlemma} that there are $2^n$ V's in $F_n$.  Therefore, there are $2^{n+1}$ degree one nodes in $F_n$ and $2^{n-1}(I_n+3)-2^{n+1}=2^{n-1}(I_n+3-2^2)=2^{n-1}(I_n-1)$ degree four nodes.  Substituting $n-1$ for $n$ in this last expression gives (\ref{crosscount}).  We note that this lemma is stated only for $n\geq2$ because $F_1$ never has a cross since there are only degree one nodes in $F_0$.  

To obtain the spectrum of $\Delta_n$ restricted to the cross, we must consider the functions in the domain of the Laplacian on the cross.  We can think of the cross as two X-shapes, (such as $F_1$ in Figure \ref{j2}) connected at their four outer nodes.  The orthogonality conditions from \cite{eigen} force the function on the bottom X to equal the negative of the function on the top X.  The value of the function on the top of the X determines the value of the function on the bottom.  The width of the X shape is $2I_n^{-1}$ and will have Dirichlet boundary conditions at the degree two nodes. Any function can be decomposed as symmetric and anti-symmetric with respect to the upper and lower branches of the X. We consider the two cases in turn.

In the symmetric case, the function is the same along the top and bottom branches of the X.  Therefore we need only to look at the top branch as it fully determines the bottom branch.  Here we are looking for eigenfunctions on an interval of length $2 I_n^{-1}$ and zero at the boundaries.  These are $\sin(\frac{1}{2} I_n k \pi x)$ with $k=1,2,\ldots$ and $x \in [0,2 I_n^{-1}]$.  The associated eigenvalues to these functions are those given in \eqref{crossspect1}.  

In the antisymmetric case, the function value horizontally along the bottom branch of the X is the negative of the value along the top branch.  At the central node, where the two branches meet, these two values must equal, so they must be zero.  We then effectively have the X broken up into two \emph{V}'s of length $I_n^{-1}$ but with Dirichlet boundary conditions at either end.  Looking at one of these \emph{V}'s, we still have the value along the bottom branch equal to the negative of the value along the top, so we consider only the top branch.  Here we look for functions of length $I_n^{-1}$ with Dirichlet boundary conditions at both ends.  This has already been done in Lemma \ref{looplemma} for the loop shape.  There we got $\sin(I_n k \pi x)$ with $k=1,2,\ldots$ and $x \in [0,I_n^{-1}]$ as the eigenfunctions and $\{ [I_n k \pi]^2 : k=1,2,\ldots \}$ as the spectrum.  This spectrum has multiplicity two because there are two halves of in the cross.    
\end{proof}

Now we must consider the graph $F_0$ and the eigenvalues it contributes to the spectrum.  This graph is just the unit interval, and has Neumann boundary conditions forced by the Kerchoff matching conditions.  So we are looking for eigenfunctions on intervals with length one and zero derivative at either end.  These come in the form $\cos(k\pi x)$ where $k=0,1,\ldots$.  This results in the following spectrum with multiplicity one:
\begin{equation}\label{linespect}
\{ [k \pi]^2 : k=0,1,\ldots \}
\end{equation}

Table \ref{spectrumtable} we summarizes the results of these lemmas.  In order to obtain the full spectrum with multiplicities, these sets must be combined with the multiplicities over all $n\geq 0$.  Hence, Theorem \ref{thm:spect} holds.
\begin{table}\center\small
\begin{tabular}{|c|c|c|c|}
\hline
Shape & Count & Spectrum & n Value \\ \hline
$F_0$ & 1 & $\{ [k \pi]^2 : k=0,1,\ldots \}$ & $n=0$ \\ \hline
V & $2^n$ & $\{ [I_n (k+\frac{1}{2}) \pi]^2:k=0,1,\ldots \}$ & $n\geq1$ \\ \hline
Loop & $2^{n-1}(j_n-2)(I_{n-1})$ & $\{ [I_n k \pi]^2 : k=1,2,\ldots \}$ & $n\geq1$  \\ \hline
\multirow{2}{*}{Cross} & \multirow{2}{*}{$2^{n-2}(I_{n-1}-1)$} & $\{ [\frac{1}{2} (I_n k \pi)]^2 : k=1,2,\ldots \}$ & \multirow{2}{*}{$n\geq2$}  \\
& & $\{ [I_n k \pi]^2 : k=1,2,\ldots \}^{\times 2}$ & \\
\hline
\end{tabular}
\caption{Summary of Lemmas \ref{vlemma} through \ref{crosslemma}}\label{spectrumtable}
\end{table}

\subsection{Numerical Computations of the Spectrum}\label{ssec:matlab}

A MATLAB script described in \cite{eigen} calculated the spectrum of the Laplacian for constant $j_n$ for all $n\in\mathbb N$ by producing the incidence matricies of the approximating graphs.  We modified this script to handle general Laakso spaces.  As in the original script, the eigenvalues are calculated using the \emph{eigs} function, which is based on ARPACK (see users guide \cite{arpack}).

Quantum graphs approximating the Laakso space associated with the sequence $j_i$=[2,3,2,3,...] is shown in Figure \ref{pic:2323} and the first twenty eigenvalues of the Laplacian on this space are shown in Table \ref{table23}. The first twenty eigenvalues of Laplacians on other Laakso spaces are shown in Table \ref{23table}.  These computations agree with the calculations found in Lemmas \ref{vlemma} through \ref{crosslemma}.

\section{Heat Kernel}\label{sec:heatmethods}

Given a Laplacian on a Laakso space, the trace of its heat kernel will be obtained following the same procedure used in \cite{physics} where the heat kernel's trace was found for the diamond fractal.  From \cite{tutorial} the heat kernel of the Laplacian is $$p(t,x,y)=\sum_{k,l,m}\psi_{k,l}(y)\psi_{k,m}(x)e^{-tE_k}$$ where $\psi_{k,l}$ and $\psi_{k,m}$ are $L^{2}-$normalized eigenfunctions of $\Delta$. The heat kernel on Laakso spaces will be further studied in \cite{futureben} where continuity and bounds will be proved. The trace of the heat kernel $Z(t)$ is defined in \cite{physics} as
\begin{equation}
Z(t)=\int p(t,x,x)dx=\sum_k g_k e^{-E_kt} \label{heatframe},
\end{equation}  
where $E_k$ are the eigenvalues of the Laplacian on the fractal and $g_k$ are the respective multiplicities associated with those eigenvalues.  Associated with the heat kernel is the \emph{spectral zeta function} also defined in \cite{physics} from the heat kernel as
\begin{equation} 
\zeta(s,\gamma)=\frac{1}{\Gamma(s)}\int_0^\infty t^{s-1} Z(t) e^{-\gamma t}dt \label{zetadefn}
\end{equation}
where $\Gamma(s)=\int_0^\infty t^{s-1}e^{-t} dt$ is the gamma function.  We set $\gamma=0$ throughout the rest of this paper.  Substitute \eqref{heatframe} into \eqref{zetadefn} to obtain 
\begin{eqnarray}
\zeta(s,\gamma)&= &\frac{1}{\Gamma(s)}\int_0^\infty t^{s-1} \displaystyle\sum_k g_k e^{-E_k t} e^{-\gamma t}dt\nonumber\\ 
&=& \displaystyle\sum_k \frac{g_k}{(E_k+\gamma )^s}. \label{zetaframe}
\end{eqnarray}
The next step in any specific example is to simplify the spectral zeta function by recognizing Riemann zeta functions, $\zeta_R(s)=\displaystyle\sum_{n=0}^\infty \frac{1}{n^s}$, and identifying the other terms as geometric series.

\begin{defn}
Define $r=\displaystyle\lim_{n\to\infty}(I_n)^{1/n}$ when this limit exists.  In the case of self-similar spaces $r$ is the contraction ratio since $I_n^{-1}$ is the diameter of each cell. There is for any value of $r$ a sequence $j_n$ that will produce that value.
\end{defn}

Once all of the series are simplfied, the poles of $\zeta(s,0)$ in the complex plane can be calculated.  The poles of $\zeta(s)$ for the diamond fractals are given in \cite{physics} as
\begin{equation}
s_m=\frac{d_h}{d_w}+\frac{2 i \pi m}{d_w \ln r},\hspace{1cm}m\in\mathbb Z,\label{diamondpole}\end{equation}
where $d_h$ and $d_s$ are the Hausdorff and walk dimensions respectively.

Since the spectral zeta function was expressed as an integral of $Z(t)$, applying an inverse Mellin transform \cite{physicshandbook} allows the heat kernel to be expressed as 
\begin{equation}
Z_D(t)=\frac{1}{2\pi i}\int_{a-i\infty}^{a+i\infty} \zeta_D(s) \Gamma(s) t^{-s}ds.\label{invmellin}\end{equation}  
By the Residue Theorem, $Z_D(t)$ is the sum of the residues of $\zeta_D(s) \Gamma(s) t^{-s}$.  The residue must also be calculated at $s=0$ (a pole for $\Gamma(s)$) and at $s=1/2$ (a pole for the $\zeta_R(2s)$ term in $\zeta(s,\gamma)$).

It is known that $\overline{\zeta_R(s)}=\zeta_R(\overline{s})$ and $\overline{\Gamma(s)}=\Gamma(\overline{s})$ for all complex $s$.  Thus, the residues from $s_m$ and $s_{-m}$ are complex conjugates of one another; therefore their sums equal twice the real part of the residue of $s_m$.  The complex values of the trace of the heat kernel yield oscillatory behavior in the heat kernel.  We shall observe what happens to the heat kernel as $t\to0$.  Therefore, only the leading term with the most negative real power of $t$ as well as any constants are included.  For example, a result of \cite{physics} shows that for the diamond fractals, the trace of the heat kernel is
\begin{equation}
Z_D(t)\sim \zeta_D(0) + \frac{r^{d_h-1}-1}{\log r^{d_w}} \frac{1}{t^{d_s/2}}(a_0 + 2 Re(a_1t^{-2 i \pi / (d_w \log r)})) +  ... .
\end{equation}

This shows that the dominating power of $t$ in the leading term as $t \to 0$ is $-d_s/2$.  This incidentally is the complex dimension introduced in \cite{Lapidus}. We shall now perform the same calculation for general Laakso spaces.

\section{The Trace of the Heat Kernel on Laakso Spaces}\label{sec:heat}
In Laakso spaces the dimensions and value of $r$ are not always known by other means.  Therefore the poles will be calculated analytically and the results will provide information about the Hausdorff, walk, and spectral dimensions.  We now give the leading powers of the trace of the heat kernel for the Laplacian on a general Laakso space.
\begin{theorem} \label{thm:heat}
For the Laakso space associated to the sequence $\{j_i\}_{i=1}^\infty$ the trace of the heat kernel is
\begin{eqnarray} \label{Laaksoheat}
Z(t)&=&\sum_{n=2}^{\infty}2^{n-1}(I_{n-1}-1)\sum_{k=1}^{\infty}e^{-k^2\pi^2I_n^2t}\notag\\&&+\sum_{n=2}^{\infty}2^{n-2}(I_{n-1}-1)\sum_{k=1}^{\infty}e^{-\frac{1}{4}k^2\pi^2I_n^2t}+\sum_{n=1}^{\infty}2^n\sum_{k=0}^{\infty}e^{-(k+1/2)^2\pi^2I_n^2t}\notag\\&&+\sum_{n=1}^{\infty}2^{n-1}I_{n-1}(j_n-2)\sum_{k=1}^{\infty}e^{-k^2\pi^2I_n^2t}+\sum_{k=0}^{\infty}e^{-k^2\pi^2t}
\end{eqnarray}
with an associated spectral zeta function
\begin{eqnarray}\label{zetawsums}
\zeta_L(s)=&\frac{\zeta_R(2s)}{\pi^{2s}}&\left[\left(\sum_{n=2}^{\infty}\frac{2^{n-1}(I_{n-1})(2^{2s-1}+j_n-1)+2^{n-1}(\frac{3}{2}2^{2s}-3)}{I_n^{2s}}\right)\right.\notag\\
&&\phantom{[[}\left.+\frac{2^{2s+1}-4+j_1}{j_1^{2s}}+1\right].
\end{eqnarray}
\end{theorem}

\begin{proof}The spectrum of the Laplacian on various Laakso spaces was given in Table \ref{spectrumtable} as
\begin{eqnarray}
\sigma(\Delta_L) =  &\displaystyle\bigcup_{n=2}^\infty \bigcup_{k=1}^\infty \left\{k^2\pi^2I_n^2\right\} \cup \bigcup_{n=2}^\infty \bigcup_{k=1}^{\infty} \left\{\frac{k^2\pi^2 I_n^2}{4}\right\}\cup \bigcup_{n=1}^\infty \bigcup_{k=0}^{\infty} \{(k+1/2)^2\pi^2 I_n^2\}\notag\\
&\cup\displaystyle\bigcup_{n=1}^\infty \bigcup_{k=1}^\infty \{k^2\pi^2I_n^2\}\cup\bigcup_{k=0}^{\infty} \pi^2 k^2
\end{eqnarray}
with respective multiplicities
\begin{equation}
2^{n-1}(I_{n-1}-1), \hspace{.25cm}2^{n-2}(I_{n-1}-1), \hspace{.25cm}2^n, \hspace{.25cm}2^{n-1}I_{n-1}(j_n-2),\hspace{.25cm}1.
\end{equation}

Direct substitution of the these values into \eqref{heatframe} gives the heat kernel.  By \eqref{zetaframe} the associated spectral zeta function is
\begin{eqnarray}
\zeta_L(s)=&\displaystyle\sum_{n=2}^{\infty}\sum_{k=1}^{\infty}\frac{2^{n-1}(I_{n-1}-1)+2^{n-2+2s}(I_{n-1}-1)}{(I_n^2k^2\pi^2)^s}\hspace{1.6cm}\notag\\
& + \displaystyle\sum_{n=1}^{\infty}\sum_{k=0}^{\infty}\frac{2^{n+2s}}{(I_n^2(2k+1)^2\pi^2)^s}+\sum_{n=1}^{\infty}\sum_{k=1}^{\infty}\frac{2^{n-1}I_{n-1}(j_n-2)}{(I_n^2k^2\pi^2)^s}\notag\\
&+\displaystyle\sum_{k=1}^{\infty}\frac{1}{(k^2\pi^2)^s}.\hspace{6.15cm}
\end{eqnarray}
This can be simplified by identifying Riemann zeta functions and, in certain cases, a Dirichlet Lambda function \cite{functionhandbook}.  Then the function can be manipulated in to one sum.

\begin{eqnarray}
\zeta_L(s)=&\frac{\zeta_R(2s)}{\pi^{2s}}&\left[\left(\sum_{n=2}^{\infty}\frac{2^{n-1}(I_{n-1})(2^{2s-1}+j_n-1)+2^{n-1}(\frac{3}{2}2^{2s}-3)}{I_n^{2s}}\right)\right.\notag\\
&&\left.+\frac{2^{2s+1}-4+j_1}{j_1^{2s}}+1\right]
\end{eqnarray}
as claimed in the proposition.
\end{proof}

This expression of the spectral zeta function cannot be simplified further for a general Laakso space with an arbitrary sequence of $j_i$'s.  However, it provides a common starting point for Laakso spaces in which the sequence of $j_i$'s is known. The next step is to locate the poles of the spectral zeta function.  Recall that only poles which yield the most negative real power of $t$ are considered since they produce the dominating behavior as $t \to 0$ in the trace of the heat kernel. There are in general more poles than these.

\begin{prop}\label{prop:poles}
Of all the poles of the spectral zeta function in the complex plane, the poles that yield the most negative real power of $t$ in the leading term of the trace of the heat kernel are located at 
\begin{equation}\label{poles}
s_m=\frac{\log2r +2\pi i m}{\log r^2}
\end{equation}
for integer $m$.
\end{prop}

\begin{proof}
The trace of the heat kernel requires the most negative power of $t$ which corresponds to the poles with the greatest real component due to the $t^{-s}$ term in \eqref{invmellin}.  The proof of the proposition relies on analyzing the series in \eqref{zetawsums} to find the poles.  Note that the in the series in \eqref{zetawsums}, the numerator has two terms; values of $s$ will be calculated that will make the denominator grow at the same rate as the numerator.  The value of $s$ with the greatest real part are the poles that will be used.  First rewrite the denominator of \eqref{zetawsums} as $I_n^{s_1}I_n^{s_2}$ where $2s=s_1+s_2$.  Select $s$ to match the rate of growth for the first term in the numerator.  To make $I_n^{s_1}$ grow at the same rate as $I_{n-1}$, $s_1$ should equal 1.  To have $I_n^{s_2}$ grow at the same rate as $2^{n-1}$ $s_2$ should equal $\log_r2$.  Therefore, for the first term, $2s=1+\log_r2$.  It can be verified that any poles from the second term in \eqref{zetawsums} will not have a real component as large as this pole.  Therefore, the real part of the poles that will yield the desired leading term in the trace of the heat kernel are $s=\frac{1}{2}+\frac{\log2}{2\log r}=\frac{\log2r}{\log r^2}$.  Including $2\pi i m$ in the numerator gives all of the complex values of this pole
\end{proof}

\begin{cor}\label{cor:sm}The dominating $t$ term in the trace of the heat kernel has power $-ds/2=-Re(s_m)$.  
\end{cor}
\begin{proof} Recall from Section \ref{sec:heatmethods} that the dominating power of $t$ in the trace of the heat kernel as $t \to 0$ was $-d_s/2$.  A result of Proposition \ref{prop:poles} is that the greatest real component of the poles yield the dominating power of $t$.  Therefore, when calculating the residue of $\zeta_L(s)\Gamma(s)t^{-s}$ at the pole obtained from \eqref{poles}, the real power of $t$ will be precisely $-Re(s_m)$.  But as stated at the beginning of the proof, it is also equal to $-d_s/2$.  Thus, $d_s/2=Re(s_m)$\end{proof}

\begin{cor}\label{cor:ds}
The spectral dimension of any Laakso space with $\{j_i\}$ such that $r$ exists is
\begin{equation} d_s=\frac{\log2r}{\log r}. \end{equation}
The walk dimension, $d_w$, is 2 for any Laakso space which implies $d_h=d_s$. \end{cor}
\begin{proof}  This follows directly from Proposition \ref{prop:poles} and Corollary \ref{cor:sm} since $d_s/2=Re(s_m)=\frac{\log2r}{\log r^2}$ which implies $d_s=\frac{\log 2r}{\log r}$.  The walk dimension is a result in \cite{futureben}.  It does agree with $d_h$ and $d_s$ via the Einstein relation $2d_h/d_w=d_s$. \end{proof}

The next two subsections give the trace of the heat kernel for two specific Laakso spaces where the sum in the spectral zeta function can be evaluated exactly: $j=2$, and $j=\{2,3,2,3...\}$.

\subsection{Laakso Space with j=2}
\begin{prop}\label{Laaksoj=2}
For the Laakso space where at each iteration $j=2$, written $L_2$, the trace of the heat kernel is

$Z_{L_2} \sim \zeta_{L_2}(0)+\frac{1}{16t \log 2}\left(1+2Re\left(  
\frac{6 \zeta_R(2+\frac{4\pi i}{\log 4})\Gamma(1+\frac{2\pi i}{\log 4})}
{t^{\frac{2\pi i}{\log 4}}\pi^{2+\frac{4\pi i}{\log 4}}}
 \right)+...\right)$.
\end{prop}

\begin{figure}[t]
\begin{center}
\includegraphics[scale=.4]{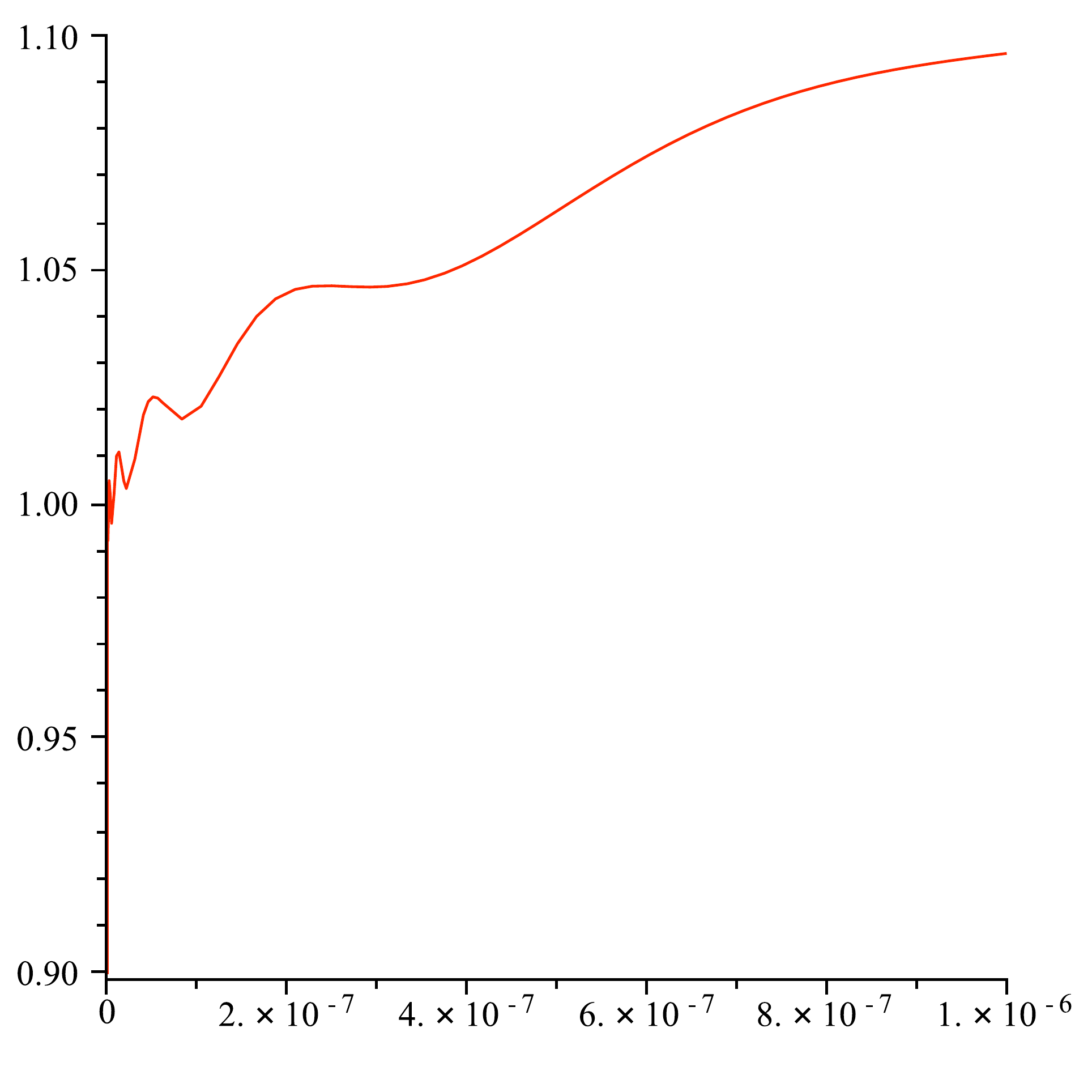}%\hspace{0cm}\includegraphics[scale=.25]{LaaksoGraph2}\\
\end{center}
\caption{Heat kernel $Z_{L_2}$, normalized by the leading non-oscillating term for the $j=2$ Laakso space. The variable $s$ is on the horizontal axis.}
\label{Graph}
\end{figure}

\begin{proof}
For Laakso spaces with a fixed $j=2$, $I_n=2^n$.  Substituting these into \eqref{Laaksoheat} and \eqref{zetawsums} gives the following two equations
\begin{eqnarray}
Z_{L_2}(t)=&\displaystyle\sum_{n=1}^\infty \left[2^{2n-2}-2^{n-1}\right] \sum_{k=1}^\infty e^{-t2^{2n}k^2\pi^2}\hspace{2.8cm}\notag\\
&+\displaystyle\sum_{n=2}^\infty \left[2^{2n-3}-2^{n-2}\right] \sum_{k=1}^\infty e^{-tk^2\pi^2 2^{2n-2}}\hspace{2.0cm}\notag\\
&+ \displaystyle\sum_{n=1}^\infty 2^n \sum_{k=0}^\infty e^{-t (2k+1)^2\pi^2 2^{2n-2}}\hspace{3.2cm}
\end{eqnarray}
and
\begin{equation}
\zeta_{L_2}(s,0)=\frac{\zeta_R(2s)}{\pi^{2s}}\left( 
\frac{4(2^{2s-1}+1)}{4^{s}(4^{2}-4)}+\frac{6(2^{s2-1}-1)}{4^{s}(4^{s}-2)} + \frac{2^{s+1}-2 + 2^{2s}}{4^{s}}
 \right)
\end{equation}
 which has poles
\begin{equation}\label{j=2poles}
z=\left( \frac{1}{2}+\frac{2 \pi i m}{\log4}\right),z=\left(1+\frac{2 \pi i m}{\log4}\right),  \forall m \in \mathbb{Z}.
\end{equation}  

By Proposition \ref{prop:poles}, only the second term of the above equation with the greatest real part contributes to the leading $t$ term.  Then an inverse Mellin transform is applied just as in Theorem \ref{thm:heat}.  Table \ref{j=2residue:} shows the residues of the integrand after the inverse Mellin transform at the poles given in \eqref{j=2poles} as well as $s=0$ from the $\Gamma(s)$ term.  Again we add complex conjugates and take only the most negative powers of t.  Notice that when adding the residue from the poles in \eqref{j=2poles} only the poles with real part one contribute to the heat kernel, the others with the exception of zero and one half have residue zero.  Once all of the residues are simplified we obtain the expression for the trace heat kernel's leading term as given in the proposition.\end{proof}

\begin{table}[hbt]
\begin{center}
  \begin{tabular}{ | c | c | }
    \hline
    $s_m$ & Residue \\ \hline
    0 & $\zeta_{L_2}(0))$ \\ \hline
    1/2 & $\frac{3}{8\sqrt{\pi}\sqrt{t}}$ \\ \hline
    1 & $\frac{1}{16t \log 2}$ \\ \hline
   $1\pm \frac{2\pi im}{\log 4}$  & $\frac{1}{16t \log 2} \frac{6\zeta_R(2+\frac{4\pi im}{\log(4)})\Gamma(1+\frac{2\pi im}{\log(4)})}{t^{\frac{2\pi im}{\log(4)}}\pi^{2+\frac{4\pi im}{\log(4)}}}$ \\ \hline
%    $1-\frac{2 i \pi}{\log 4}$  & $\frac{1}{8t \log 2} t^{2 i \pi / (\log 4)} \Gamma(1-\frac{2 i \pi}{\log 4})$ \\ \hline
    $\frac{1}{2}+\frac{2 i \pi}{\log 4}$  & 0 for $m \neq 0$ \\ \hline
%    $\frac{1}{2}-\frac{2 i \pi}{\log 4}$  & $\frac{1}{8 \sqrt{t} \log 2} t^{2 i \pi / (\log 4)} a_{-1}$ \\ \hline
  \end{tabular}
\end{center}
\caption{Residues of the integrand of the inverse Mellin transfrom for given poles of the spectral zeta function for Laakso spaces with a fixed $j=2$}
\label{j=2residue:}
\end{table}

The power of 1 of t in the denominator of the leading term in the proposition implies that the spectral dimension $d_s$ for $L_2$ is 2.  Knowing that the Hausdorff dimension $d_h=2$ from Table \ref{hdim:}, we conclude that the walk dimension $d_w=2$ since $d_s=2d_h/d_w$.
Notice the poles and dimensions of this fractal were explicitly calculated.  But Corollary \ref{cor:ds} and Proposition \ref{prop:poles} yield the same result once $r=\displaystyle \lim_{n \to \infty}I_n^{1/n}=\lim_{n \to \infty}(2^n)^{1/n} =2$ is known.

\subsection{j=\{2,3,2,3...\} Laakso Space}
\begin{prop}
For the Laakso space with $j_{2k}=3$ and $j_{2k-1}=2$ where $k\ge 1,$ the trace of the heat kernel is 
{\small
\begin{align}
Z_L(t) \sim & \zeta_L(0) +\\ \notag &\frac{1}{24t^{\frac{1}{2}+\frac{\log(2)}{\log(6)}}\log(6)}\left(\frac{\Gamma\left(\frac{1}{2} + \frac{\log(2)}{\log(6)}\right) \zeta_R\left( 1 + \frac{2\log(2)}{\log(6)} \right) }{\pi^{1 +\frac{2\log(2)}{\log(6)}} }+\right. \\ \notag &\left. \sum_{m=1}^{\infty}2Re\left( \frac{\Gamma\left( \frac{1}{2}+\frac{\log(2)}{\log(6)}+\frac{2\pi im}{\log(6)} \right) \zeta_R\left( 1 +\frac{2\log(2)}{\log(6)} + \frac{4\pi im}{\log(6)} \right)}{\pi^{1 +\frac{2\log(2)}{\log(6)} + \frac{4\pi im}{\log(6)}} t^{-\frac{2\pi i m}{\log(6)}} }\right)\right. \\ \notag & \left.\left( \frac{2^{2+\frac{4\log(2)}{\log(6)}+\frac{8\pi im}{\log(6)}} + 10 \cdot 2^{1+\frac{2\log(2)}{\log(6)}+\frac{4\pi im}{\log(6)}} + 12}{2^{1+\frac{2\log(2)}{\log(6)} +\frac{4\pi im}{\log(6)}}} \right)  \right).
\end{align}
}
\end{prop}

\begin{proof}
In this case $r=\displaystyle\lim_{n\to\infty}I_n^{1/n}=\displaystyle\lim_{n\to\infty}(2^{n/2}3^{n/2})^{1/n}=\sqrt6$.  This locates the poles with largest real part at $s_m=\frac{1}{2}+\frac{\log(2)}{\log6}+\frac{2\pi i m}{\log6}$.

The next step is to obtain and simplify the spectral zeta function associated with the trace of heat kernel given in \eqref{zetawsums}.  Since $j_i$ alternates between 2 and 3, the following values  can be directly substituted in, for any $k$ we have:
\begin{center}
  $\begin{array}{  c  c  c  c  c }
  I_1=2&I_2=6&I_{2k-2}=6^{k-1}&I_{2k-1}=2\times6^{k-1}&I_{2k}=6^k
 \end{array}$
\end{center}

In preparation for substituting these values into \eqref{zetawsums} the sum is split into two sums, one over even $n$ and the other over odd 

{\small
\begin{align}
& \zeta_L(s)=\notag \\
\sum_{n=2}^{\infty}&\frac{2^{n-1}(I_{n-1}-1)+2^{n-2+2s}(I_{n-1}-1)+2^{n+2s}-2^n+2^{n-1}I_{n-1}(j_n-2)}{I_n^{2s}}\notag\\
=&\frac{2(I_1-1)+2^{2s}(I_1-1)+2^{2+2s}-2^2+2I_1}{I_2^{2s}}\notag\\
\phantom{=}&+\displaystyle\sum_{k=2}^\infty \bigg(\frac{2^{2k-2}(I_{2k-2}-1)+2^{2k-3+2s}(I_{2k-2}-1)+2^{2k-1+2s}-2^{2k-1}}{I_{2k-1}^{2s}}\notag\\
\phantom{=}&\phantom{}+\frac{2^{2k}(I_2k-1)+2^{2k-2+2s}(I_{2k-1}-1)+2^{2k+2s}-2^{2k}+2^{2k-1}(I_{n-1})}{I_{2k}^{2s}}\bigg).
\end{align}  
}

Substituting the known values  of $I_{2k},\ I_{2k+1},\ j_{2k},$ and $j_{2k+1}$ we obtain
\begin{align}
\zeta_L(s)=\frac{\zeta_R(2s)}{\pi^{2s}}&\left[ \left(\frac{2^{2s}+4}{12} + \frac{s^{2s-1}+1}{2^{2s}} \right)\frac{24}{6^{2s}-24} +  \right. \notag \\
& \left. \left( 3\cdot 2^{2s}-6 + \frac{6 \cdot 2^{2s}-12}{2^{2s}} \right)\frac{1}{6^{2s}-4} \right].
\end{align}

To apply the inverse Mellin transform as was done in \eqref{invmellin} to obtain an expression for $Z_L(t)$.  Then $Z_L(t)$ can be calculated using the sum of the residues of $\zeta_L(s)\Gamma(s)t^{-s}$ using the poles obtained in \eqref{poles} as well as $s=1/2$ (obtained from the Riemann zeta function, $\zeta_R(2s)$ and $s=0$ (from $\Gamma(s)$).  Table \ref{23residue:} lists the residues for those poles.

\begin{table}[hbt] \small
\begin{center}
  \begin{tabular}{ | c | c | } 
    \hline
    $s_m$ & Residue \\ \hline
    $0$&$\zeta_L(0)$\\ \hline
    $\frac{1}{2}$&$\frac{2}{3\sqrt{\pi t}}$\\ \hline
    $s_m = \frac{1}{2}+\frac{\log(2)}{\log6}+\frac{2\pi im}{\log(6)}$&$\frac{\Gamma(s_m)\zeta_R(2s_m)}{t^{s_m}\pi^{2s_m}} \frac{2^{-2s_m}}{24\log(6)}\left( 2^{4s_m} + 10 \cdot 2^{2s_m} + 12 \right) $\\ \hline
   $s_k = \frac{\log(2)}{\log(6)} + \frac{2\pi i m}{\log(6)}$&$\frac{\Gamma(s_k)\zeta_R(2s_k)}{t^{s_k}\pi^{2s_k}}\frac{3}{8\log(6)} \frac{4^{2s_k}-4}{4^{s_k}}$\\ \hline
  \end{tabular}
\end{center}  
\caption{Residues of the integrand of the inverse Mellin transfrom for given poles of the spectral zeta function for the \{2,3,2,3...\}Laakso spaces}
\label{23residue:}
\end{table}

The sum of the residues in Table \ref{23residue:} expresses the trace of the heat kernel with the leading terms as shown in the statement of the proposition.  Note that the general terms in $Z_L(t)$ are shown in the proposition to indicate the behavior of the trace of the heat kernel in more detail.
\end{proof}

\begin{cor}
The exponent of $t$, $\frac{\log(2\sqrt6)}{\log{6}}$, in the leading term as $t$ goes towards zero implies that the spectral dimension of this Laakso space is $d_s=\log24/\log6$.  The Hausdorff dimension for this fractal is given in Table \ref{hdim:} as $d_h=\log24/\log6$.  \end{cor}
Again, the same results are obtained by simply applying Corollary \ref{cor:ds} and Proposition \ref{prop:poles}.

\section*{Acknowledgments}
The authors would like to thank Alexander Teplyaev, Luke Rogers, Robert Strichartz, Shotaro Makisumi, Grace Stadnyk, Jun Kigami, and Naotaka Kajino.

\end{document}